\documentclass[reqno,10pt]{amsart}

\usepackage{amssymb,amsmath,amsthm,amsfonts}
\usepackage{mathrsfs,dsfont, comment,mathscinet}
\usepackage{enumerate,esint}
\usepackage[colorlinks=true, linkcolor=blue, citecolor=blue, urlcolor=blue]{hyperref}
\usepackage[left=2.8cm,right=2.8cm,top=2.8cm,bottom=2.8cm]{geometry}
\usepackage{mathtools}
\usepackage[pdftex]{color, graphicx}
\usepackage{bbm}
\usepackage[english]{babel}
\usepackage{comment} 
\usepackage{esint}
\usepackage{pdfsync}
\usepackage{enumitem}
\usepackage{caption} 
\usepackage{epstopdf}
\allowdisplaybreaks

\numberwithin{equation}{section}
\theoremstyle{plain}

\newtheorem{thm}{Theorem}[section]
\newtheorem{prop}[thm]{Proposition}
\newtheorem{lemma}[thm]{Lemma}

\newtheorem{definition}[thm]{Definition}
\theoremstyle{definition}
\newtheorem{remark}[thm]{Remark}

\newcommand{\N}{\mathbb{N}}

\newcommand{\Z}{\mathbb{Z}}

\newcommand{\R}{\mathbb{R}}
\newcommand{\e}{\varepsilon}

\newcommand{\dx}{\,\mathrm{d}x}
\newcommand{\dy}{\,\mathrm{d}y}
\newcommand{\dif}{\,\mathrm{d}}

\def\hom{{\rm hom}}
\newcommand{\per}{{\rm per}}
\newcommand{\umin}{u^{\rm min}}

\title{First-order homogenization}
\author{Riccardo Cristoferi}
\address{Department of Mathematics, IMAPP, Radboud University, Heyendaalseweg, 6525 AJ 
Nijmegen, The Netherlands}
\email{cristoferi@science.ru.nl}

\author{Lorenza D'Elia}
\address{Institute of Analysis and Scientific Computing, TU Wien, Wiedner Hauptstraße 8-10, 1040
Vienna, Austria}
\email{lorenza.delia@tuwien.ac.at}
\date{}

\begin{document}

\begin{abstract}
We provide a first-order homogenization result for quadratic functionals.
In particular, we identify the scaling of the energy and the explicit form of the limiting functional in terms of the first-order correctors.
The main novelty of the paper is the use of the dual correspondence between quadratic functionals and PDEs, combined with a refinement of the classical Riemann-Lebesgue Lemma.
\end{abstract}

\maketitle


\section{Introduction}

First-order homogenization does not exist. The non-existence has its roots in two types of boundary effects. The first one arises when the domain is not a disjoint union of suitable rescaled copies of the periodicity cell (see \cite[Example 1.12]{BrTr08}). The second comes from the oscillatory nature of correctors in the homogenization theory (see \cite[Equation (2.12)]{AllAmar}). Nevertheless, in this manuscript, we will provide a first-order $\Gamma$-convergence result for quadratic energies under suitable assumptions which allow us to \textquoteleft forget\textquoteright\; about the boundary. In particular, the scaling and the principal part of the energy in the bulk will be identified. 
\smallskip

Nowadays, homogenization is a very well-established mathematical theory describing how the microstructure affects the overall behavior of a material (see, e.g., \cite{BR18, BD98}).
The mathematical literature on the topic is too vast for an exhaustive list, and thus we limit ourselves to mention here some examples where it has been successfully applied: from thin structures (see, for instance, \cite{BalAlOmn, BraDEl, BouFra, ZhiPya}), to phase separation (see, for instance, \cite{CriFonGan-fixed_super, CriFonGan_dep, CriFonHagPop}), and from micromagnetism (see \cite{AloDiF, ChrKre, DavDiF, DavDElIng}) to supremal functionals (see \cite{BriPriGar, DElEleZap}).

The prototypical example is a family of functionals $\mathcal F_\e: L^p(\Omega; \R^M)\to \R\cup\{+\infty\}$, for $p\in(1,\infty)$, with $\e>0$ being the length scale characterizing the fine structure, of the form
               \begin{equation}
                   \notag
                   \mathcal F_\e(u)\coloneqq
                   \begin{dcases}
                       \int_\Omega W\left({x\over \e}, \nabla u(x)\right)\dx & \mbox{if } u\in W^{1,p}_0(\Omega; \R^M),\\
                       +\infty & \mbox{else}.
                   \end{dcases}
               \end{equation}
The first variable of the energy density $W$ accounts for the presence of a periodic microstructure, which is reflected in requiring that $W(\cdot, \xi)$ is a periodic function.  
The variational investigation of the periodic homogenization goes back to the end of Seventies. In \cite{Marcellini}, the limiting functional $\mathcal F_{\hom}$ of $\mathcal F_\e$ has been fully characterized in the scalar case, i.e. $M=1$, and under the assumptions of convexity and of $p$-growth of $W(x, \cdot)$. The $\Gamma$-limit $F_{\hom}: L^p(\Omega)\to \R\cup\{+\infty\}$ is given by 
                 \begin{equation}
                   \label{def:limitfunperhom}
                   \mathcal F_{\hom}(u)\coloneqq
                   \begin{dcases}
                       \int_\Omega W_{\hom}\left(\nabla u(x)\right)\dx & \mbox{if } u\in W^{1,p}_0(\Omega),\\
                       +\infty & \mbox{else},
                   \end{dcases}
               \end{equation}  
with the effective energy density $W_{\hom}:\R^N\to \R$ being characterized through the so-called {\it cell formula}
      \begin{equation}
          \label{def:cellformula}
          W_{\hom}(\xi)\coloneqq \inf\left\{ \int_{[0,1)^N}  W(y, \xi+\nabla u(x)): u\in W_0^{1,p}([0,1)^N)   \right\}.
      \end{equation}
Removing the assumption of convexity and in the vectorial framework, the analysis of the $\Gamma$-limit has been carried out independently in \cite{Br85} and \cite{Mueller}. In this case, the limiting energy  $\mathcal F_{\hom}: L^p(\Omega; \R^M)\to \R\cup\{+\infty\}$ is again of the form \eqref{def:limitfunperhom} but the homogenized energy density  $W_{\hom}:\R^{M\times N}\to \R$ is characterized by the {\it asymptotic cell formula}
        \begin{equation}
          \label{def:asyformula}
          W_{\hom}(\Xi)\coloneqq \lim_{k\to\infty}{1\over k^N}\inf\left\{ \int_{[0,k)^N}  W(y, \Xi+\nabla u(x)): u\in W_0^{1,p}([0,k)^N; \R^M)   \right\}.
      \end{equation}
It is worth noticing that in the scalar setting and assuming the convexity of $W$ in the second variable, formula  \eqref{def:asyformula} turns into \eqref{def:cellformula} (see \cite[Lemma 4.1]{Mueller}).
\smallskip

The aim of the present paper is to undertake a first-order analysis of a suitable version functional $\mathcal F_\e$ via $\Gamma$-convergence.
We focus on quadratic energies $F_\e: L^2(\Omega)\to \R\cup\{+\infty\}$ of the form 
      \begin{equation}
          \notag
          F_\e(u)\coloneqq 
              \begin{dcases}
                  \int_\Omega A\left({x\over \e}\right)\nabla u(x)\cdot \nabla u(x)\dx - \int_\Omega f(x)u(x)\dx &\mbox{if } u\in H^1_0(\Omega),\\
                  +\infty & \mbox{else},
              \end{dcases}
      \end{equation}
where $\Omega\subset\R^N$ is a bounded open set, $f\in L^2(\Omega)$ and $A:\Omega\to\R^{N\times N}$ is a matrix-valued function in $L^\infty$ that is $[0,1)^N$-periodic, symmetric and with lower and upper quadratic bounds.
To identify the contribution of the bulk in the first-order limit, we mimic the approach deployed by Allaire and Amar in \cite{AllAmar} assuming that 
   \begin{equation*}
       \Omega \coloneqq [0,1)^N\;\;\;\;\; \hbox{and}\;\;\;\;\; \e= {1\over n},\;\hbox{with } n\in \N\setminus\{0\}. 
   \end{equation*}
This allows us to get rid of the first type of boundary effect. To tackle the second issue with boundary effect, we restrict the admissible class for the source term $f$ (see Section \ref{sec:mainresult} for further details). 
Using the asymptotic expansion of functionals given in \cite{AB93} in terms of $\Gamma$-convergence, we consider the functionals $F_{\e}^1: L^2(\Omega)\to \R\cup\{+\infty\}$ defined as
        \begin{equation}
            \notag
             F_{\e}^1 (u)\coloneqq {F_{\e}(u) - \min_{H^1_0(\Omega)} F_{\hom} \over \e}.
        \end{equation}
The principal result of the present manuscript is the identification of the scale $\e$ above as well as the $\Gamma$-limit of $F_{\e}^1$ (see Theorem \ref{thm:mainresult}). The main novelty lies in the use of the dual correspondence between quadratic functionals and PDEs. To the best of the authors' knowledge, this is the first time that these two theories are combined together to get a variational result. 
\smallskip

We briefly outline the strategy we employ. The unique minimizer of $F_{\e}$ turns out to be the unique solution of the following elliptic problem with Dirichlet boundary conditions
            \begin{equation*}
                \left\{
                \begin{aligned}
                    -{\rm div}\left(A\left({x\over\e}\right)\nabla u_\e(x)\right)&=f(x) && \mbox{in  } \Omega,\\
                    u_\e(x)&=0 && \mbox{on } \partial\Omega.
                \end{aligned}
                \right.
            \end{equation*} 
The investigation of such an elliptic problem has been broadly carried out by many authors, \cite{Allaire, BLP78, JKO94} to name a few (see \cite{BR18} for an extensive review on the topic). To get a homogenized equation, the classical strategy relies on the two-scale expansion developed in \cite{Allaire, N89} of the solution $u_\e$ : 
        \begin{equation}
        \label{twoscaleexp}
            u_\e(x)= u_0(x) + \e u_1\left(x, {x\over\e}\right) + \dots,
        \end{equation}
where $u_0$ is the solution of the homogenized equation given by 
   \begin{equation*}
                \left\{
                \begin{aligned}
                     -{\rm div}\left(A_\hom\nabla u_0(x)\right)&=f(x) && \mbox{in  } \Omega,\\
                     u_0(x)&=0 && \mbox{on } \partial\Omega, 
                \end{aligned}
                \right.
            \end{equation*} 
where $A_{\hom}$ is defined through the cell formula \eqref{def:cellformula} with $W(y, \xi+\nabla u)= A(y)(\xi+\nabla u(x))\cdot(\xi+\nabla u(x))$.
Moreover, the function $u_1$ is defined through the first-order correctors (see Section \ref{sect:preliminaries} for the precise definition). 
One would be tempted to use the ansatz \eqref{twoscaleexp}  in the variational analysis for the functional $F_\e^1$ to deduce the limiting energy. However, this idea does not work out. The reason is that the following estimate
             \begin{equation*}
                 \left\| u_\e - u_0 -  \e u_1\left(\cdot, {\cdot\over\e}\right) \right\|_{H^1(\Omega)}\leq C\sqrt{\e}
             \end{equation*}
turns out to be sharp (see, e.g., \cite{BLP78}). This surprising result suggests the presence of another phenomenon, known as {\it boundary layers}. They are further first-order corrections needed to match the boundary conditions (see \cite{AllAmar, AKMetal, GM10}). Due to the high oscillatory nature of these functions, their energy contribution is not clearly quantifiable with respect to the parameter $\e$.
This is what we have referred to as the second type of boundary effect.
In order to avoid this high oscillatory behavior at the boundary, we essentially consider the case where the function $u_\e$ is compactly supported in $\Omega$ by requiring the source term $f$ to be in a specific class (see Assumption {\rm (H5)} in Section \ref{sec:mainresult}).
This enables us to get the first-order $\Gamma$-limit by using the ansatz in \eqref{twoscaleexp} together with a refinement of the classical Riemann-Lebesgue Lemma (see Proposition \ref{prop:RL}).
\smallskip

Finally, we show that the first-order $\Gamma$-convergence analysis is not needed when the functional $F_\e$ only depends on the function $u$ and not on its gradient  $\nabla u$ (see Theorem \ref{thm:main_Lp}). Indeed, in such a case, we have that the value of the minimum of the functional $F_\e$ is the same as the one of $F_\hom$.
This implies that the expansion by $\Gamma$-convergence does not provide additional information on the minimizers of the functional $F_\e$.
\smallskip

The paper is organized as follows. In Section \ref{sect:setup} we specify the set-up of the problem. The preliminaries and the technical results are given in Section \ref{sect:preliminaries} and Section \ref{sect:techlemma}, respectively
We then turn to the proofs of the main result: in Section \ref{sect:compactness} we provide the compactness, while Section \ref{sect:liminf} and Section \ref{sect:limsup} are devoted to the lower and upper bound, respectively.
Finally, Section \ref{sect:generalfunc} is devoted to the proof of the first-order homogenization for functionals in $L^p$, for $p\in(1,\infty)$.


\section{Set-up of the problem and main result}
\label{sect:setup}

\subsection{Basic notation} Here we collect the basic notation we are going to use throughout the manuscript.  Let $Y\subset\R^N$ be a periodicity cell, namely 
    \[
    Y = \left\{\, \sum_{i=1}^N \lambda_i v_i \,:\, 0<\lambda_i<1,\, \sum_{i=1}^N \lambda_i=1 \,\right\},
    \]
where $v_1,\dots,v_N$ is a basis of $\R^N$.
Without loss of generality, up to a translation, we can even assume that $Y$ has its baricenter at the origin. This assumption is just to simplify the writing of the main result. \\
The space $H^1_{\per}(Y)$ is the subset of $H^1(Y)$ of functions with periodic boundary conditions.
More precisely, $u\in H^1(Y)$ if and only if the function $\widetilde{u}: \R^N\to\R$ defined as $\widetilde{u}(y)\coloneqq u(\widetilde{y})$ belongs to $H^1_{\mathrm{loc}}(\R^N)$, where
\[
y=\sum_{i=1}^N \lambda_i v_i,\quad\quad\quad
\widetilde{y}\coloneqq\sum_{i=1}^N \{\lambda_i\} v_i,
\]
and $\{\lambda_i\}\coloneqq \lambda_i - \lfloor \lambda_i \rfloor$.\\
Given a function $f$, the notation $f^\e$ stands for $f^\e(x)\coloneqq f(x/\e)$.
Moreover, we denote by $\partial_i$ the $i^{th}$ partial derivative operator with respect to the variable $x$, and by $\partial_{y_i}$ the $i^{th}$ partial derivative operator with respect to the variable $y$. In particular, we have that
\[
\partial_i f^\e(x) = \frac{1}{\e} \partial_{y_i}f^\e(x).
\]
Finally, the symbol $\langle\cdot\rangle_Y$ denotes the average over $Y$, i.e.,
\[
\langle f\rangle_Y\coloneqq {1\over |Y|}\int_Y f(y)\dif y,
\]
with $|Y|$ being the $N$-dimensional Lebesgue measure of $Y$.


\subsection{Main result}\label{sec:mainresult}
Let $\Omega\subset\R^N$ be an open, bounded set, and 
let $f\in L^2(\Omega)$.
Let $A:\R^N\to\R^{N\times N}$ be a matrix-valued function in $L^\infty$ such that
    \begin{itemize}
        \item [{\rm (H1)}] $A$ is $Y$-periodic;
        \item [{\rm (H2)}] $A$ is symmetric, i.e., $a_{ij}(y)=a_{ji}(y)$;
        \item [{\rm (H3)}] there exist two positive constants $\alpha, \beta$ such that 
             \begin{equation}
                 \notag
                 \alpha|\xi|^2\leq A(y)\xi\cdot\xi\leq \beta |\xi|^2,
             \end{equation}
             for all $\xi\in\R^N$.
    \end{itemize}
For $\e>0$, let $F_{\e}: L^2(\Omega)\to \R\cup \{+\infty\}$ be the functional defined as
     \begin{equation}
         \notag
         F_{\e}(u)\coloneqq \int_\Omega A^{\e}\left(x\right) \nabla u(x)\cdot\nabla u(x) \dx
            - \int_\Omega f(x) u(x)\dx,
     \end{equation}
if $u\in H^1_0(\Omega)$, and as $F_{\e}(u)\coloneqq+\infty$ otherwise in $L^2(\Omega)$.

Under Assumptions {\rm (H1)-(H3)}, we know (see, for instance, \cite{Marcellini}, \cite[Theorem 1.3]{Mueller}, or \cite[Corollary 24.5]{DM}) that $\{F_{\e}\}_\e$ $\Gamma$-converges with respect to the weak topology of $H^1(\Omega)$, or equivalently the strong topology of $L^2(\Omega)$, to the effective functional $F^0_\hom: L^2(\Omega)\to\R\cup\{+\infty\}$ given by
     \begin{equation}
         \notag
         F^0_\hom(u)\coloneqq \int_{\Omega} A_{\hom}\nabla u(x)\cdot\nabla u(x) \dx
          - \int_\Omega f(x) u(x)\dx,
     \end{equation}
if $u\in H^1_0(\Omega)$, and by $F^0_\hom(u)\coloneqq+\infty$ otherwise in $L^2(\Omega)$.
Here, the effective matrix $A_{\hom}$ is a constant matrix given by the {\it cell-formula}
    \begin{align}
        \label{eq:def_A_hom}
        A_{\hom}\xi\cdot\xi \coloneqq \inf\left\{ \int_{Y} A(y)   (\xi + \nabla \phi(y))\cdot(\xi+ \nabla \phi(y)) \dif y : \phi \in H^1_{\per}(Y) \right\}.
    \end{align}
We refer to this $\Gamma$-convergence result as the zeroth-order term in the expansion by $\Gamma$-convergence of $F_\e$.
For our analysis, we need to recall the following. Using the fact that the functional in \eqref{eq:def_A_hom} is quadratic, for each $\xi\in\R^N$, the minimization problem defining $A_{\hom}\xi\cdot\xi$ has a unique solution (up to an additive constant), denoted by $\psi_\xi$.

    \begin{definition}\label{def:first_order_correctors}
    Let $e_1,\dots,e_N$ be the standard orthonormal basis of $\R^N$.
    For each $i=1,\dots,N$, let $\psi_i\in H^1_{\per}(Y)$ be the unique solution to
    \[
    \inf\left\{ \int_{Y} A(y)   (e_i + \nabla \phi(y))\cdot(e_i+ \nabla \phi(y)) \dif y : \phi \in H^1_{\per}(Y),\,\, \int_Y \varphi(y)\dif y =0 \right\}.
    \]
    The function $\psi_i$ is called the \emph{first-order corrector for $A$} associated to the vector $e_i$.
    \end{definition}

    \begin{remark}
    It turns out that the map $\xi\mapsto \psi_\xi$ is linear (see \cite[Example 25.5]{DM}).
    Namely, 
        \[
        \psi_\xi = \sum_{i=1}^N \psi_i \xi_i,
        \]
    where $\xi=(\xi_1,\dots,\xi_N)$. Therefore, the knowledge of the first-order correctors for $A$ is sufficient to obtain $A_{\hom}$.
    \end{remark}

The goal of this paper is to develop further the expansion by $\Gamma$-convergence of $F_\e$  (see \cite{AB93} for further details).
As explained in the Introduction, there are two issues with boundary effects in obtaining such a result.
Thus, in order to carry out  our analysis, we need to assume the following
   \begin{itemize}
       \item[{\rm (H4)}] $\Omega \coloneqq Y$ and we take the sequence $\e_n= {1\over n}$, with $n\in \N\setminus\{0\}$;
       \item[{\rm (H5)}] It holds that
       \[
       f = -\mathrm{div}(A_{\hom}\nabla g),
       \]
       for some $g\in C^\infty_c(Y)$.
   \end{itemize}
Assumption {\rm (H4)} is to ensure that $\Omega$ can be obtained as a disjoint union of copies of rescaled versions of the periodicity cell $Y$.
On the other hand, assumption {\rm (H5)} is in order to avoid using boundary layers, which have two main issues: they do not have a variational definition, and their contribution to the energy is not clearly quantifiable in terms of the parameter $\e$, which prevents us from getting an order of the energy.
From now on, recalling assumption {\rm (H4)}, we will denote by $F_n$ the functional $F_{\e_n}$.

We are now in position to write the asymptotic expansion through $\Gamma$-convergence we will study.

    \begin{definition}
    For $n\in\N\setminus\{0\}$, we define the functional $F_n^1: L^2(\Omega)\to \R\cup\{+\infty\}$ as
         \begin{equation}
             \label{eq:F_e_1}
             F_n^1 (u)\coloneqq {F_{n}(u) - \min_{H^1_0(\Omega)}F^0_{\hom} \over \e_n}.
         \end{equation}
    \end{definition}

\begin{remark}
Note that, using standard estimates (see the proof of Proposition \ref{prop:cpt}), it is possible to prove that, for each $n\in\mathbb N\setminus\{0\}$, the minimization problem
    \[
    \min_{u\in H^1_0(\Omega)} F_n(u)
    \]
admits a unique solution $\umin_n\in H^1(\Omega)$. 
In a similar way, it is possible to prove that the minimization problem
    \[
    \min_{u\in H^1_0(\Omega)} F^0_\hom(u)
    \]
admits a unique minimizer $\umin_0\in H^1(\Omega)$.
 In particular, we have that
 \[
 F_{n}^1 (u) = {F_{n}(u) - F^0_{\hom}(\umin_0) \over \e_n}.
 \]
Moreover, the $\Gamma$-convergence of $F_{n}$ to $F^0_\hom$ together with the compactness of $\{\umin_n\}_{n}$ yields that $\umin_n\to\umin_0$ in $L^2(\Omega)$ as $n\to\infty$.
\end{remark}

\begin{remark}\label{rem:u0min}
   Note that assumption {\rm (H5)} implies that $u^{\min}_0=g$. In particular, $u^{\min}_0\in C^\infty_c(\Omega)$.
\end{remark}

Our goal is to study the $\Gamma$-limit of the family of functionals $\{F_n^1\}_{n}$.
To this end, we introduce the candidate limiting functional.
    \begin{definition} 
    Define the functional $F^1_{\hom}: L^2(\Omega)\to \R\cup\{+\infty\}$ as
          \begin{align}
               \label{eq:F_hom_1}
              F^1_{\hom}(u)&\coloneqq
              \sum_{i,j=1}^N \int_\Omega
    \nabla(\partial_i \umin_0\partial_j \umin_0)(x)\dx
    \cdot\left\langle y \left[ a_{ij}
            + 2A e_j \cdot \nabla \psi_i
            + A \nabla \psi_i \cdot \nabla \psi_j \right]\right\rangle_Y \nonumber \\
            &\hspace{1cm}+2 \sum_{i=1}^N \int_\Omega   \langle \psi_i A\rangle_Y \nabla \umin_0(x)\cdot \partial_i\umin_0(x)\dx \nonumber \\
            &\hspace{1.5cm}+ 2 \sum_{j=1}^N\sum_{i=1}^N \int_{\Omega} \partial_j \umin_0(x) \langle \psi_i A\nabla_y\psi_j \rangle_Y  \cdot \partial_i \nabla \umin_0(x) \dx \nonumber \\
            &\hspace{2cm}-\sum_{j=1}^N \langle \psi_j \rangle_Y \int_\Omega f(x) \partial_j u_0^{\min}(x)\dx
           \end{align}
    if $u=u^{\min}_0$, and $+\infty$ else in $L^2(\Omega)$.
    Here, $\psi_i$ are the first-order correctors defined in  Definition \ref{def:first_order_correctors}.
    \end{definition}

We now state the main result of the present paper.
     
    \begin{thm}
    \label{thm:mainresult}
    Let $F_n^1$ and $F_{\hom}^1$ be the functionals given by \eqref{eq:F_e_1} and \eqref{eq:F_hom_1}, respectively.
    Assume that Assumptions {\rm (H1)}-{\rm (H5)} hold.
    Then, we have the following:
       \begin{itemize}
           \item[{\rm (i)}]  Let $\{u_n\}_{n}$ be a sequence in $H^1(\Omega)$ such that $\sup_{n} F^1_n(u_n)<\infty$.
           Then, $\{u_n\}_n$ converges in $L^2(\Omega)$ to $\umin_0$, where $\umin_0$ is the unique minimizer of $F_{\hom}^0$.
           \item[{\rm (ii)}] The sequence $\{F_n^1\}_{n}$ $\Gamma$-converges with respect to $L^2(\Omega)$ topology to $F_{\hom}^1$.
       \end{itemize}
    \end{thm}

It is worth noticing that the above result can be read as follows: in the bulk, the order of the energy $F_n$ is $\e_n$, and the normalized difference from $\min F^0_{\hom}$ is given by the (constant) functional $F^1_{\hom}$.

We show that the analogous first-order $\Gamma$-expansion is trivial in the case of functionals defined on $L^p$.
We are able to prove this statement in a more general setting than that considered above.
Fix $p\in(1,\infty)$ and $M\geq1$. Let $V:\Omega\times\R^M\to\R$ be a Carath\'{e}odory function such that
    \begin{itemize}
    \item[{\rm (A1)}] For each $z\in\R^M$, the function $x\mapsto V(x,z)$ is $Y$-periodic;
    \item[{\rm (A2)}] There exists $0<c_1<c_2<+\infty$ such that
    \[
    c_1(|z|^p-1) \leq V(x,z)\leq c_2(|z|^p+1),
    \]
    for all $x\in\Omega$ and $z\in\R^M$.
    \end{itemize}

    Note that in this case the source term would only be a part of the function $V$. That is why there is no analogue of assumption {\rm (H5)}.

    \begin{definition}
    For $n\in\N\setminus\{0\}$, define $G_n: L^p(\Omega;\R^M)\to\R\cup\{+\infty\}$ as
    \[
    G_n(u) \coloneqq \int_\Omega V\left( \frac{x}{\varepsilon_n}, u(x) \right)\dx.
    \]    
    \end{definition}

    \begin{remark}
    Note that there is no loss of generality in assuming the function $V$ to be convex in the second variable. Indeed, the relaxation of $G_n$ with respect to the weak-$L^p(\Omega)$ topology is given by
    \[
    \overline{G_n}(u) = \int_\Omega V^c\left( \frac{x}{\varepsilon_n}, u(x) \right)\dx
    \]
    where $V^c$ denotes the convex envelope of the function $V$ in the second variable.
    \end{remark}

    We now introduce the homogenized functional.

    \begin{definition}
    For $z\in\R^M$ , let
    \[
    V_\hom(z)\coloneqq \inf\left\{\frac{1}{|Y|} \int_Y V^c(y,z + \varphi(y))\dif y \,: \, \varphi\in L^p(Y;\R^M),\, \int_Y \varphi(y)\dif y = 0 \right\}.
    \]
    Define the functional $G_\hom: L^p(\Omega;\R^M)\to\R\cup\{+\infty\}$ as
    \[
    G_\hom \coloneqq \int_\Omega V_\hom (u(x))\dx.
    \]    
    \end{definition}

    Using the same arguments as in the proof of \cite[Theorem 3.3]{CriFonGan_dep}, we get the following.

    \begin{lemma}
	Assume that Assumptions {\rm (A1)}-{\rm (A2)} hold.
    Then, $G_n\to G_{\hom}$ with respect to the weak-$L^p(\Omega)$ topology.
    \end{lemma}

Next result justifies our claim that the $\Gamma$-expansion is trivial in such a case.

    \begin{thm}
    \label{thm:main_Lp}
    Assume that Assumptions {\rm (A1)}-{\rm (A2)}, as well as {\rm{(H4)}} hold and that, for all $y\in Y$, the function $u\mapsto V(y,u)$ is strictly convex.
	Let $m\in\R^M$.
    Then,
    \[
    \min\left\{ G_n(v) : v\in L^p(\Omega;\R^M),\, \int_\Omega u(x)\dx =m \right\}
    =
    \min\left\{ G_{\hom}(v) : v\in L^p(\Omega;\R^M),\, \int_\Omega u(x)\dx =m \right\},
    \]
    for all $n\in\N\setminus\{0\}$.
    \end{thm}
    
    \begin{remark}
    The assumption of strict convexity is needed only to simplify the strategy of the proof, since we need to invert the relation
    \[
    \partial_u V(x,v)=c.
    \]
    Strict convexity gives us a unique inverse, while with convexity alone we would have to get a slightly more involved argument. Since this result is to show that for these functionals there is no need to a first order $\Gamma$-expansion, we decided to use this extra technical assumption.
    \end{remark}
    
    \begin{remark}
    Note that, in this case, boundary conditions are not natural. This is why we consider a mass constraint instead.
    \end{remark}


\section{Preliminaries}
\label{sect:preliminaries}

\subsection{Recall of homogenization}
We recall the foundations of the homogenization theory for elliptic equations. Even if this theory is classical by now, we revisit it since we will use a slightly different definition of the correctors than that in  \cite{AllAmar, BLP78}. \\
Let $\Omega$ be a bounded and open subset of $\R^N$. Let $A$ be the matrix-valued function satisfying Assumptions {\rm (H1)-(H3)}. For a given $f\in L^2(\Omega)$, we consider the following equation
       \begin{equation}
       	\label{eq:diveq}
        \left\{
       	\begin{aligned}
       		-{\rm div}(A^\e(x)\nabla u_\e(x))
            &= f(x) && \mbox{in }\Omega,\\
       		u_\e(x)&=0 && \mbox{on } \partial\Omega.
       	\end{aligned}
        \right.
       \end{equation}
It is well-known that this problem admits a unique solution in $H^1_0(\Omega)$ (see, e.g., \cite{BLP78}). To carry out a homogenization procedure, the solution $u_\e$ is assumed to have the following two-scale expansion
    \[
    	u_\e(x) = u_0\left(x, {x\over \e}\right) + \e u_1\left(x, {x\over \e}\right) + \e^2u_2\left(x, {x\over \e}\right) + \dots,  
    \]
where each $u_i$ is $Y$-periodic with respect to the fast variable $y={x\over\e}$. The variables $x$ and $y$ are treated as independent. Plugging such an expansion into \eqref{eq:diveq} and identifying powers of $\e$, we get a cascade of equations. Here, we only care about the second-order expansion.  Therefore, setting $\mathcal{A}_\e \phi(x)\coloneqq -{\rm div}(A^\e(x)\nabla \phi(x))$, we deduce that 
\[
\mathcal{A}_\e= \e^{-2}\mathcal{A}_0+\e^{-1}\mathcal{A}_1+ \mathcal{A}_2,
\]
where
     \begin{align}
     	\mathcal{A}_0 &\coloneqq -\sum_{i=1}^N{\partial_ {y_i}}\left(\sum_{j=1}^N a_{ij}(y){\partial_{y_j}}\right),\notag\\
     	\mathcal{A}_1 &\coloneqq -\sum_{i=1}^N{\partial_{y_i}}\left(\sum_{j=1}^N a_{ij}(y){\partial_j}\right) - \sum_{i=1}^N{\partial_i}\left(\sum_{j=1}^N a_{ij}(y){\partial_{y_j}}\right),\notag\\
     	\mathcal{A}_2 &\coloneqq -\sum_{i=1}^N{\partial_i}\left(\sum_{j=1}^N a_{ij}(y){\partial_j}\right)\notag.
     \end{align} 
Using \eqref{eq:diveq}, matching powers of $\e$ up to second order gives us the following equations
       \begin{align}
       	\mathcal{A}_0u_0&=0,\label{firsteq}\\
       	\mathcal{A}_0u_1 + \mathcal{A}_1u_0&=0,\label{secondeq}\\
       	\mathcal{A}_0u_2 + \mathcal{A}_1u_1 + \mathcal{A}_2u_0&=f,\label{thirdeq}\\
            \mathcal{A}_0u_3 + \mathcal{A}_1u_3 + \mathcal{A}_2u_1&=0.\label{fortheq}
       \end{align}
From \eqref{firsteq}, it follows that $u_0(x,y)\equiv u_0(x)$.
The solution $u_1$ to \eqref{secondeq} is given by 
      \begin{equation}
      	\label{defu1}
      	u_1\left(x, {x\over \e}\right)= \sum_{j=1}^N \psi_j^{\e}\left(x\right) {\partial_j}u_0(x) + \widetilde{u}_1(x),
      \end{equation}
 where, for $j=1,\dots, N$, $\psi_j$ is the unique solution in $H^1_{\rm per}(Y)$ to the problem
     \begin{equation}
     	\label{eq:first_order_corrector}
     	\begin{dcases}
     		\mathcal{A}_0 \psi_j(y)= \sum_{i=1}^N\partial_{y_i} a_{ij}(y) & \mbox{in } Y,\\
     		\int_Y \psi_j(y)\dif y=0, \\
     		y\mapsto \psi_j(y) &Y\mbox{-periodic}.
     	\end{dcases}
     \end{equation}
Namely, for $j=1,\dots, N$, the function $\psi_j$ are the first-order corrector defined in Definition \ref{def:first_order_correctors}. The solution to equation \eqref{thirdeq} is given by 
     \begin{equation}
     	\label{def:u2}
     	u_2\left(x, {x\over \e}  \right) = \sum_{i,j=1}^N\chi^{\e}_{ij}\left(x\right) {\partial^2_{ij}} u_0(x)
        +\sum_{j=1}^N \psi_j^{\e}\left(x\right) {\partial_j}\widetilde{u}_1(x)+ \widetilde{u}_2(x),
     \end{equation}
where for all $i,j=1,\dots, n$, the {\it second-order corrector} $\chi_{ij}\in H^1_{\per}(Y)$ satisfies the following auxiliary problem
    \begin{equation}
    	\label{def:second_order_corrector}
    	\begin{dcases}
    		\mathcal{A}_0 \chi_{ij}(y) = b_{ij}(y)-\int_Y b_{ij}(y)\dif y &\mbox{in } Y,\\
    		\int_Y \chi_{ij}(y)\dif y=0,\\
             y\mapsto \chi_{ij}(y)&Y\mbox{-periodic},
    	\end{dcases}
    \end{equation}
with   
    \begin{equation}
    	\label{eq:defbij}
    	b_{ij}(y)\coloneqq a_{ij}(y) +\sum_{k=1}^N a_{ik}(y){\partial_{y_k} \psi_j}(y) + \sum_{k=1}^N{\partial_{y_k}}(a_{ki}\psi_j)(y).
    \end{equation}
It is worth recalling that the function $\widetilde{u}_1$ in \eqref{defu1} can be taken identically equal to zero if one is only interested in the first-order expansion of $u_\e$. Otherwise, $\widetilde{u}_1$ is determined by the compatibility condition of equation \eqref{fortheq}, namely,
        \begin{equation}
            \label{def:problemtildeu1}
            {\rm div}(A_{\hom}\nabla \widetilde{u}_1(x))= -\sum_{i,j,k=1}^N c_{ijk}{\partial^3_{ijk} u_0^{\min}}(x),
        \end{equation}
where $A_{\hom}$ is the homogenized matrix defined by \eqref{eq:def_A_hom} and for $i,j,k=1,\dots,N$, the constant $c_{ijk}$ is defined as
        \begin{equation}
            \label{def:cijk}
            c_{ijk}\coloneqq \left\langle  \sum_{l=1}^N a_{kl}{\partial_{y_l} \chi_{ij}} + a_{ij}\psi_k \right\rangle_Y. 
        \end{equation}
Since in the present paper, we are interested in expanding $u_\e$ up to the second-order, the function $\widetilde{u}_2$ in \eqref{def:u2} can be taken identically equal to zero.

\begin{remark}
In the following, we will use the above theory for the function $u_0^{\min}$ in place of $u_0$.
\end{remark}


\subsection{Periodic functions}
We collect some useful properties of the space $L^2_{\per}(Y; \R^N)$ of the periodic functions (see \cite[Chapter 1]{JKO94} for further details).\\
The space $L^2_{\rm sol}(Y)$ of solenoidal periodic functions is defined as
   \begin{equation*}
       L^2_{\rm sol}(Y)\coloneqq\{f\in L^2_{\per}(Y; \R^N)\; : \; \hbox{div}f=0\hbox{ in } \R^N \},
   \end{equation*}
where the equality $\hbox{div}f=0$ is to be intended in distributional sense, i.e., 
   \begin{equation*}
       \int_{\R^N} f(x)\cdot \nabla \varphi(x)\; dx=0, \;\;\;\hbox{for any } \varphi \in C^\infty_c(\R^N).
   \end{equation*}
The space $L^2_{\rm sol}(Y)$ turns out to be a closed subspace of $L^2(Y; \R^N)$. Setting 
    \begin{equation}
        \notag
        \mathcal V^2_{\rm pot}(Y)\coloneqq \{\nabla u\;: u\in H^1_{\per}(Y; \R^N)\},
    \end{equation}
we immediately deduce  the following orthogonal representation
         \begin{equation}
             \notag
             L^2_{\per}(Y; \R^N) = L^2_{\rm sol}(Y) \oplus \mathcal V^2_{\rm pot}(Y). 
         \end{equation}
The next proposition provides a representation of solenoidal functions, which will be a key tool in Proposition \ref{prop:estimate}. 

\begin{prop}
\label{prop:skewsymmatrix}
    Let $f\in L^2_{\rm sol}(Y)$. Then, $f=(f_1, \dots, f_N)$ can be represented in the form
              \begin{equation}
                  \notag
                  f_j(x)= \langle f_j\rangle_Y +\sum_{i=1}^N {\partial_i} \alpha_{ij}(x),\;\;\;\hbox{for all } j=1,\dots, N,
              \end{equation}
    where $\alpha_{ij}\in H^1(Y)$ is such that $\alpha_{ij}=-\alpha_{ji}$ and $\langle \alpha_{ij}\rangle_Y=0$, for all $i,j=1, \dots, N$.      
   
\end{prop}

  \begin{proof}
     Without loss of generality, we can assume that $Y=[0, 2\pi)^N$. Using the Fourier series of $f$, we have that
           \begin{equation}
               \notag
               f=\langle f \rangle_Y + \sum_{\substack{k\in\Z^N\\k\neq 0}} f^k e^{\mathrm{i}k\cdot x},
           \end{equation}
    where $f^k$ is the Fourier coefficient of $f$. We claim that $f^k\cdot k=0$, for each $k\in\Z^N$. Indeed, for fixed $k\in \Z^N\setminus\{0\}$,  we can decompose $f^k$ as  $f^k = C_1^k k + C_2^k v$, for some constants $C_1^k, C_2^k$ and $v\in ({\rm span}(k))^{\perp}$. Therefore, 
        \begin{equation}
            \notag
            f=\langle f\rangle_Y + \sum_{\substack{k\in\Z^N\\k\neq 0}}C_1^k ke^{\mathrm{i}k\cdot x} + \sum_{\substack{k\in\Z^N\\k\neq 0}}C_2^k v e^{\mathrm{i}k\cdot x}.
        \end{equation}
       Noticing that $\mathrm{i} ke^{\mathrm{i}k\cdot x}=\nabla(e^{\mathrm{i}k\cdot x}) \in \mathcal V^2_{\rm pot}(Y)$ and since by assumption, we know that $f\in L^2_{\rm sol}(Y)$, we obtain that $C_1^k=0$, for any $k\in\Z^N\setminus\{0\}$. This leads us to conclude that $f^k\cdot k=0$, and that $C_2^k v = f^k$. \\
       Therefore,
       \[
       f^k = f^k - \frac{f^k\cdot k}{|k|^2}k,
       \]
       which writes, component-wise, as
            \begin{equation*}
                f_j^k = \sum_{i=1}^N\ g_{ij}^k k_i, \;\;\;\hbox{with }\,\,\,\,\, g_{ij}\coloneqq{f_j^k k_i- f_i^k k_j \over |k|^2}.
            \end{equation*}
        Therefore, for each $j=1,\dots, N$, we get that
        \begin{align*}
        f_j &= \langle f_j \rangle_Y + \sum_{\substack{k\in\Z^N\\k\neq 0}} f^k_j e^{\mathrm{i}k\cdot x}
        = \langle f_j \rangle_Y + \sum_{\substack{k\in\Z^N\\k\neq 0}} \sum_{i=1}^N\ g_{ij}^k k_i e^{\mathrm{i}k\cdot x}
        = \langle f_j \rangle_Y + \sum_{i=1}^N \sum_{\substack{k\in\Z^N\\k\neq 0}} \ g_{ij}^k k_i e^{\mathrm{i}k\cdot x}.
        \end{align*}
        Thus, by defining
              \begin{equation}
                  \notag
                  \alpha_{ij}(x) \coloneqq -\mathrm{i} \sum_{\substack{k\in\Z^N\\k\neq 0}} g_{ij}^k e^{\mathrm{i}k\cdot x},
              \end{equation}
        we get the desired result.
  \end{proof}


\section{Technical lemmata}
\label{sect:techlemma}

\subsection{First order Riemann-Lebesgue Lemma}

We prove a quantitative version of the Riemann-Lebesgue Lemma.
Since we will use this latter result several times, we recall it in here for the reader's convenience.

\begin{lemma}[Riemann-Lebesgue Lemma]\label{lem:RL}
Let $p\in[1,\infty)$, and let $\Omega\subset\R^N$ be an open bounded set.
Let $g\in L^p(\R^N)$ be an $Y$-periodic function.
Then, 
\[
\lim_{n\to\infty}\int_\Omega g(nx)\varphi(x)\, dx
    = \langle g \rangle_Y \int_\Omega \varphi(x)\, dx,
\]
for all $\varphi\in L^{p'}(\Omega)$, with $\tfrac{1}{p}+\tfrac{1}{p'}=1$.
\end{lemma}

\begin{remark}
The above result holds also for $p=+\infty$, but in that case we need the test function $\varphi\in L^1(\Omega)$.
\end{remark}

The refined result requires more stringent assumptions. Nevertheless, we will later show how to relax some of them.

\begin{prop}[First order Riemann-Lebesgue Lemma]\label{prop:RL}
Let $p\in[1,\infty)$, and let $g\in L^p(\R^N)$ be an $Y$-periodic function.
Let $\Omega\subset\R^N$ be an open bounded set such that, up to a set of Lebesgue measure zero, can be written as
\begin{equation}\label{eq:assumption_Omega_RL}
\Omega = \bigcup_{i=1}^k (x_i + Y),\quad x_i\in \Z^N,\quad\quad
    (x_i+Y)\cap(x_j+Y)=\emptyset\, \text{ if }\, i\neq j.
\end{equation}
Then,
\begin{align*}
&\lim_{n\to\infty} n\left[ \int_\Omega g(nx)\varphi(x)\, dx
    - \langle g \rangle_Y \int_\Omega \varphi(x)\, dx \right] = \int_\Omega \nabla \varphi(x) \dx \cdot
        \left[ \langle y g \rangle_Y - \langle y \rangle_Y\, \langle g\rangle_Y  \right],
\end{align*}
for all $\varphi\in C^2(\Omega)\cap C(\overline{\Omega})$.
\end{prop}

\begin{proof}
Let $\{z_i\}_{i}$ be an enumeration of $\Z^N$.
For each $n\in\N\setminus\{0\}$, let 
\[
I_n\coloneqq \left\{ i\in\N : \frac{1}{n} (z_i + Y)\subset \Omega \right\}.
\]
Note that, by \eqref{eq:assumption_Omega_RL}, we get that (up to a set of Lebesgue measure zero)
\[
\Omega = \bigcup_{i\in I_n} \frac{1}{n}(z_i + Y)
\]
being a disjoint union. Moreover, the cardinality of $I_n$ is $n^N|\Omega|/|Y|$.
Recall that, by periodicity of $g$, we have that $g(z_i + y) = g(y)$ for all $z_i\in\Z^N$.
Then,
\begin{align*}
\int_\Omega & g(nx)\varphi(x)\, dx
    - \frac{1}{|Y|}\int_Y g(y)\, dy \, \int_\Omega \varphi(x)\, dx  \\
& =\sum_{i\in I_n} \frac{1}{n^N }\left[ \int_Y g(z_i+y)\varphi\left(\frac{z_i + y}{n}\right)\, dy
    -  \frac{1}{|Y|}\int_Y \varphi\left(\frac{z_i + y}{n}\right)\, dy \, \int_Y g(y)\, dy  \right] \\
&=\sum_{i\in I_n} \frac{1}{n^N} \left[ \int_Y g(y)
    \left[\varphi\left(\frac{z_i + y}{n}\right) - \varphi\left(\frac{z_i}{n}\right) \right]\, dy
    -  \frac{1}{|Y|}\int_Y \left[\varphi\left(\frac{z_i + y}{n}\right) - \varphi\left(\frac{z_i}{n}\right) \right]\, dy \, \int_Y g(y)\, dy  \right].
\end{align*}
We now use a second-order Taylor expansion for $\varphi$ to write
\[
\varphi\left(\frac{z_i + y}{n}\right) - \varphi\left(\frac{z_i}{n}\right)
    = \frac{1}{n} \nabla \varphi\left(\frac{z_i}{n}\right)\cdot y 
        + \frac{1}{2n^2} H\varphi\left(\frac{z_i}{n}\right)[y,y]
        + o\left(\frac{1}{n^2}\right),
\]
where $H\varphi(x)$ denotes the Hessian matrix of $\varphi$ at $x$.
Thus, we get that
\begin{align*}
\int_\Omega g(nx)\varphi(x)\, dx
    &- \frac{1}{|Y|}\int_Y g(y)\, dy \, \int_Y \varphi(x)\, dx  \\
& = \frac{1}{n} \sum_{i\in I_n} \frac{1}{n^N} \nabla \varphi\left(\frac{z_i}{n}\right)
    \cdot \left[ \int_Y g(y)y\, dy - \frac{1}{|Y|}\int_Y y\, dy\, \int_Y g(y)\, dy  \right]
    +O\left(\frac{1}{n^2}\right).
\end{align*}
Sending $n\to\infty$, and noting that 
\[
\sum_{i\in I_n} |Y| \frac{1}{n^N}  \nabla \varphi\left(\frac{z_i}{n}\right)
\]
is the Riemann sum for the integral of $\nabla\varphi$ over $\Omega$, gives the desired result.
\end{proof}

\begin{remark}
Note that the above result is invariant under translation of the periodicity grid and of the function $g$. In particular, there is no loss of generality in assuming that $\langle y\rangle_Y=0$, namely that the barycenter of the periodicity cell is at the origin.
\end{remark}

\begin{remark}
In particular, the above result implies that
\[
\lim_{n\to\infty} n\left[ \int_\Omega g(nx)\varphi(x)\, dx
    - \langle g\rangle_Y \, \int_\Omega \varphi(x)\, dx \right] =0,
\]
for each $Y$-periodic function $\varphi\in C^2(\R^N)$.
\end{remark}

Finally, we present some extensions of the above result. Their proofs follow the same argument used above, and therefore we will not repeat them in here.

The first extension uses the fact that, in dimension $N\geq 3$, given a general open set $\Omega\subset\R^N$ with Lipschitz boundary (or even whose boundary has finite $\mathcal{H}^{N-1}$ measure), the contribution of the cubes $( z_i + Y)/n$ that intersect the boundary of $\Omega$ is of order $n^{1-N}\ll n^{-1}$.

\begin{prop}
Let $N\geq 3$, and let $\Omega\subset\R^N$ be an open set with Lipschitz boundary.
Let $p\in[1,\infty)$, and let $g\in L^p(\R^N)$ be an $Y$-periodic function.
Then,
\begin{align*}
&\lim_{n\to\infty} n\left[ \int_\Omega g(nx)\varphi(x)\, dx
    - \langle g \rangle_Y \int_\Omega \varphi(x)\, dx \right] = \int_\Omega \nabla \varphi(x) \dx \cdot
        \left[ \langle y g \rangle_Y - \langle y \rangle_Y\, \langle g\rangle_Y  \right],
\end{align*}
for all $\varphi\in C^2(\Omega)\cap C(\overline{\Omega})$.
\end{prop}

\begin{remark}
Note that, in dimension $N=2$, the contribution of the terms at the boundary of $\Omega$ is of the order $n^{-1}$. It is therefore, in general, not possible to estimate this term.
\end{remark}

Finally, we can extend the above result to the nonlinear case.

\begin{prop}
Let $p\in[1,\infty)$.
Let $V:\R^N\times\R\to\R$ be a Carath\'{e}odory function that is $Y$-periodic in the second variable, and such that $t\mapsto V(x,t)$ is of class $C^2$ for all $x\in\R^N$.
Assume that \eqref{eq:assumption_Omega_RL} hold.
Then,
\begin{align*}
&\lim_{n\to\infty} n\left[ \int_\Omega V(nx,\varphi(x))\, dx
    - \int_\Omega <V(\cdot, \varphi(x)>_Y \, dx \right] \\
&\hspace{2cm}= \int_\Omega
    \frac{1}{|Y|} \int_Y  \partial_t V(y, \varphi(x)) \nabla \varphi(x)\cdot y \dy\dx 
    - \frac{1}{|Y|^2} \int_\Omega\int_Y \int_Y 
            \partial_t V(z,\varphi(x))\nabla\varphi(x)\cdot y  \, \mathrm{d}z \dy \dx,
\end{align*}
for all $\varphi\in C^2(\Omega)\cap C(\overline{\Omega})$.
\end{prop}


\subsection{Estimates}
In this section we prove the fundamental estimate that allows us to consider the `surrogate' sequence given by the expansion with first and second-order correctors in place of the minimizer $u^{min}_\e$.
Recall that $u_0^{\min}\in C^\infty_c(\Omega)$ (see Remark \ref{rem:u0min}).

    \begin{prop}
    \label{prop:estimate}
    Let $u^{(2)}_\e$ be the function defined as
         \begin{equation}
             \label{def:ue2}
             u^{(2)}_\e(x, y) \coloneqq u_0^{\min}(x)+ \e u_1\left(x, y\right)+ \e^2 u_2\left(x, y\right),
         \end{equation}
    with $u_1$ and $u_2$ being defined as in \eqref{defu1} and \eqref{def:u2}, respectively. Then, it holds that
           \begin{equation}
               \notag
               \|{\rm div}(A^\e \nabla u^{(2)}_\e)-{\rm div}(A^{\hom} \nabla \umin_0)  \|_{H^{-1}(\Omega)}\leq C \e^2,
           \end{equation}
           for some $C<+\infty$.
    \end{prop}
    \begin{proof} First, note that 
         \begin{equation}
             \notag
             {\rm div}(A^\e(x) \nabla u^{(2)}(x)_\e-A^{\hom} \nabla \umin_0(x)) = \sum_{i=1}^N \partial_i\left( (A^\e(x)\nabla u^{(2)}_\e(x)- A^{\hom}\nabla\umin_0(x))_i\right),
         \end{equation}
        where, for $i=1,\dots, N$, the $i$-th component $(A^\e(x)\nabla u^{(2)}_\e(x)- A^{\hom}\nabla\umin_0(x))_i$ is given by 
          \begin{align}
          \label{eq:componentith}
              (A^\e(x)\nabla u^{(2)}_\e(x)- A^{\hom}\nabla\umin_0(x))_i= P_{i, \e}(x)+ \e Q_{i, \e}(x)+ \e^2R_{i,\e}(x),
          \end{align}
        with 
              \begin{align}
                  P_{i, \e}(x) &\coloneqq \sum_{j=1}^N a^\e_{ij}(x) \partial_j \umin_0(x) + \sum_{j,k=1}^N a_{ik}^\e(x)\partial_{y_k}\psi^\e_j(x)\partial_j \umin_0 (x)-
                  \sum_{j=1}^N a_{ij}^{\hom}\partial_j \umin_0(x),\notag\\
                  Q_{i, \e}(x)&\coloneqq \sum_{j,k=1}^N \psi^\e_j(x) a^\e_{ik}(x)\partial^2_{kj} \umin_0(x) + \sum_{k,j,l=1}^N a_{il}^\e(x) \partial_{y_l}\chi^{\e}_{kj}(x)\partial^2_{kj}\umin_0(x)  + \sum_{j=1}^N a_{ij}^\e(x)\partial_j\widetilde{u}_1(x) \notag\\
                  &\quad + \sum_{j,k=1}^N a_{ik}^\e(x)\partial_{y_k}\psi^\e_j(x)\partial_j\widetilde{u}_1(x),\label{eq:defQie}\\
                  R_{i, \e}(x)&\coloneqq \sum_{k,j,l=1}^N \chi_{kj}^\e(x) a_{il}^\e(x)\partial^3_{lkj}\umin_0(x) + \sum_{j,k=1}^N \psi_j^\e(x) a^{\e}_{ik}(x)\partial^2_{kj}\widetilde{u}_1(x) \notag.
              \end{align}
      We now estimate the above terms, starting from $P_{i, \e}(x)$. Note that $P_{i, \e}(x)$ rewrites as 
           \begin{equation}
               \notag
               P_{i, \e}(x) = \sum_{j=1}^N g_i^{j, \e}(x) \partial_j \umin_0(x), \;\;\; \mbox{for all } i=1,\dots, N,
           \end{equation}
       where $g_i^j$ is defined as 
            \begin{align}
                \notag
                g_i^j(y) \coloneqq a_{ij}(y)+ \sum_{k=1}^Na_{ik}(y)\partial_{y_k}\psi_{j}(y) -a_{ij}^{\hom}, \;\;\; \hbox{for all } i, j=1,\cdots, N.
            \end{align}
        For fixed $j=1,\dots, N$, set $G^j\coloneqq(g_i^1,\dots, g_i^N)$. We have that $G^j\in L^2_{\rm sol}(Y)$. Indeed, thanks to problem \eqref{eq:first_order_corrector} satisfied by $\psi_j$, it follows that 
              \begin{equation}
                  \notag
                  \mbox{div } G^j(y)= \sum_{i=1}^N \partial_{y_i} g_{i}^j(y)=0.
              \end{equation}
        Therefore, by applying Proposition \ref{prop:skewsymmatrix}, the components of $G^j$ are represented by  
               \begin{equation}
                   \label{eq:representationgij}
                   g_i^j(y) = \sum_{k=1}^N \partial_{y_k}\alpha^j_{ik}(y).
               \end{equation}
        Note that $\langle g_i^j\rangle_Y=0$, since
              \begin{equation}
                  \label{eq:defahom}
                 a_{ij}^{\hom}= \left\langle a_{ij}(y)+ \sum_{k=1}^Na_{ik}(y)\partial_{y_k}\psi_{j}(y) \right\rangle_Y. 
              \end{equation}
        Using the representation \eqref{eq:representationgij} as well as the Leibniz rule, we get that   
                 \begin{align}
                     P_{i, \e}(x) &= \sum_{j,k=1}^N \partial_{y_k} \alpha_{ik}^{j, \e}(x)\partial_j\umin_0(x) \notag\\
                     &= \e\left(\sum_{j,k=1}^N \partial_{k}(\alpha_{ik}^{j,\e}\partial_j\umin_0)(x) - \sum_{j,k=1}^N \alpha_{ik}^{j, \e}(x)\partial^2_{kj}\umin_0(x) \right).\notag
                 \end{align} 
   Hence, \eqref{eq:componentith} turns into
          \begin{align}
          \label{eq:componentith-1}
              (A^\e(x)\nabla u^{(2)}_\e(x)- A^{\hom}\nabla\umin_0(x))_i= \e\left( \sum_{j,k=1}^N \partial_{k}(\alpha_{ik}^{j,\e}\partial_j\umin_0)(x) + \widetilde{Q}_{i,\e}(x)\right)+ \e^2R_{i,\e}(x),  
          \end{align}
   where
       \begin{equation}
           \notag
           \widetilde{Q}_{i,\e}(x) \coloneqq Q_{i,\e}(x)- \sum_{j,k=1}^N \alpha_{ik}^{j, \e}(x)\partial^2_{kj}\umin_0(x),
       \end{equation}
   with $Q_{i,\e}$ being defined as in \eqref{eq:defQie}. We now estimate $\widetilde{Q}_{i,\e}(x)$. As for $P_{i,\e}$, $\widetilde{Q}_{i,\e}(x)$ can be rewritten as 
             \begin{align}
                \widetilde{Q}_{i,\e}(x)= \sum_{j,k=1}^N h_i^{jk, \e}(x)\partial_{kj}^2\umin_0(x) +\sum_{j=1}^N t_i^{j,\e}(x)\partial_j\widetilde{u}_1(x), \;\;\; \hbox{for all } i=1,\dots, N, 
             \end{align}
      where for any fixed $j,k=1,\dots, N$, the functions $h^{kj}_i$ and $t^j_i$ are defined as
            \begin{equation}
                \label{def:hikj}
                h_i^{kj} (y)\coloneqq   \psi_j(y) a_{ik}(y) + \sum_{l=1}^N a_{il}(y) \partial_{y_l}\chi_{kj}(y) -\alpha_{ik}^j(y),
            \end{equation}
            \begin{equation}
                \label{def:tij}
                t_i^j(y)\coloneqq a_{ij}(y) +\sum_{k=1}^N a_{ik}(y) \partial_{y_k}\psi_j(y).
            \end{equation}
        We claim that $H^{kj}\coloneqq (h_1^{kj}, \dots, h_N^{kj})$ as well as $T^j\coloneqq (t_1^j, \dots, t_N^j)$ belong to $L^2_{\rm sol}(Y)$. Indeed, bearing in mind that  $\chi_{kj}^\e$ is the solution to \eqref{def:second_order_corrector} together with the fact that $$\sum_{i=1}^N\partial_{y_i}(-\alpha_{ik}^j)=\sum_{i=1}^N\partial_{y_i}(\alpha_{ki}^j) = g_k^j,$$ we deduce that for fixed $j,k=1,\dots, N$,
        \begin{align}
            {\rm div}(H^{kj}) &= \sum_{i=1}^N \partial_{y_i} h_i^{kj}(y)= \sum_{i=1}^N\partial_{y_i}(\psi_j^\e a_{ik}^\e)(y) +\sum_{i=1}^N\partial_{y_i}\left(\sum_{l=1}^N a_{il}\partial_l\chi_{kj}\right)(y) +\sum_{i=1}^N \partial_{y_i}(-\alpha_{ik}^j)(y) \notag\\
            &=\sum_{i=1}^N \partial_{y_i}(\psi_j(y) a_{ik}(y)) -b_{kj}(y)+\int_Y b_{kj}dy +   g_k^j(y) =0,  \notag
        \end{align}
    where we have used the equality $a_{kj}^{\hom}=\langle b_{kj}\rangle_Y$ (cf. \eqref{eq:defbij} and \eqref{eq:defahom}). 
    Likewise, since $\psi_j$ is the solution to \eqref{eq:first_order_corrector}, we immediately conclude that 
          \begin{equation}
              \notag
              {\rm div}(T^j) = \sum_{i=1}^N \partial_{y_i} t_i^j (y)=\sum_{i=1}^N\partial_{y_i}\left(a_{ij}(y) +\sum_{k=1}^N a_{ik}(y)\partial_{y_k}\psi_j(y)\right) =0.
          \end{equation}
    Therefore, applying Proposition \ref{prop:skewsymmatrix}, it follows that the components of $H^{kj}$ and $T^j$ are represented by
             \begin{align}\notag
                 h_i^{kj}(y)= \langle h_i^{kj}\rangle_Y+ \sum_{l=1}^N \partial_{y_l} \beta_{il}^{kj}(y)\;\;\;\hbox{and}\;\;\;
                  t_i^j(y)= \langle t_i^{j}\rangle_Y+\sum_{k=1}^N \partial_{y_k}\gamma_{ik}^j(y), 
             \end{align}
    for all $i=1,\dots, N$. This implies that 
           \begin{align}
             \widetilde{Q}_{i,\e}(x) &= \sum_{j,k=1}^N \langle h_i^{kj}\rangle_Y\partial^2_{kj}u_0^{\min}(x)+ \sum_{j,k,l=1}^N \partial_{y_l} \beta_{il}^{kj,\e}(x)\partial^2_{kj}u_0^{\min}(x) \notag\\
             &\quad+ \sum_{j=1}^N\langle t_i^{j}\rangle_Y\partial_j\widetilde{u}_1(x)+\sum_{j,k=1}^N \sum_{k=1}^N \partial_{y_k}\gamma_{ik}^{j, \e}(x) 
  \partial_j\widetilde{u}_1(x)\notag\\
  & =\sum_{j,k=1}^N \langle h_i^{kj}\rangle_Y\partial^2_{kj}u_0^{\min}(x)+ \sum_{j=1}^N\langle t_i^{j}\rangle_Y\partial_j\widetilde{u}_1(x) \notag\\
  &\quad +\e\biggl( \sum_{j,k,l=1}^N \partial_l(\beta_{il}^{kj, \e}\partial_{kj}^2\umin_0)(x)- \sum_{j,k,l=1}^N \beta_{il}^{kj,\e}(x)\partial_{lkj}^3\umin_0(x)\notag\\
  &\quad+ \sum_{j,k=1}^N\partial_i(\gamma_{ik}^{j,\e}\partial_j\widetilde{u}_1)(x) -\sum_{j,k=1}^N\gamma_{ik}^{j,\e}(x)\partial^2_{kj}\widetilde{u}_1)(x)      \biggr)\notag
           \end{align}
From \eqref{eq:componentith} together with \eqref{eq:componentith-1}, we conclude that 
     \begin{align}
         (A^\e(x)\nabla u^{(2)}_\e(x)- A^{\hom}\nabla\umin_0(x))_i&= \e\left( \sum_{j,k=1}^N \partial_{k}(\alpha_{ik}^{j,\e}\partial_j\umin_0)(x)+\sum_{j,k=1}^N \langle h_i^{kj}\rangle_Y\partial^2_{kj}u_0^{\min}(x)+ \sum_{j=1}^N\langle t_i^{j}\rangle_Y\partial_j\widetilde{u}_1(x)\right)\notag\\
         &\quad + \e^2\left( \sum_{j,k,l=1}^N \partial_l(\beta_{il}^{kj, \e}\partial_{kj}^2\umin_0)(x)+ \sum_{j,k=1}^N\partial_i(\gamma_{ik}^{j,\e}\partial_j\widetilde{u}_1)(x) +  \widetilde{R}_{i,\e}(x) \right),\notag
     \end{align}
    where $\widetilde{R}_{i,\e}(x)$ is defined as 
         \begin{align}
             \widetilde{R}_{i,\e}(x)&\coloneqq R_{i,\e}(x) - \sum_{j,k,l=1}^N \beta_{il}^{kj,\e}(x)\partial_{lkj}^3\umin_0(x)-\sum_{j,k=1}^N\gamma_{ik}^{j,\e}(x)\partial^2_{kj}\widetilde{u}_1(x). \notag
         \end{align}
    Noticing that $\langle h_i^{kj}\rangle_Y= c_{ijk}$ (cf. \eqref{def:cijk} and  \eqref{def:hikj}) and $\langle t_i^{j}\rangle_Y= a_{ij}^{\hom}$ (cf. \eqref{eq:defahom} and \eqref{def:tij}) as well as bearing in mind problem \eqref{def:problemtildeu1}, we obtain that
          \begin{align}
              \notag
              \sum_{i=1}^N \partial_i\left(\sum_{j,k=1}^N \langle h_i^{kj}\rangle_Y\partial^2_{kj}u_0^{\min}(x)+ \sum_{j=1}^N\langle t_i^{j}\rangle_Y\partial_j\widetilde{u}_1(x) \right)=0.
          \end{align}
    Moreover, since $\alpha_{ik}^j=-\alpha_{ki}^j$, it follows that 
         \begin{equation}
             \notag
            \sum_{i=1}^N\partial_i\left(\sum_{j,k=1}^N \partial_{k}(\alpha_{ik}^{j,\e}\partial_j\umin_0)(x)
             \right)=0.             
         \end{equation}
    Likewise, since $\beta_{il}^{kj} = - \beta_{li}^{kj}$ and $\gamma_{ik}^{j}=-\gamma_{ki}^{j}$,
          \begin{equation}
              \notag
              \sum_{i=1}^N\partial_i\left(
              \sum_{j,k,l=1}^N \partial_l(\beta_{il}^{kj, \e}\partial_{kj}^2\umin_0)(x)\right)= \sum_{i=1}^N\partial_i\left(\sum_{j,k=1}^N\partial_i(\gamma_{ik}^{j,\e}\partial_j\widetilde{u}_1)(x)\right)=0.              
          \end{equation}
    Therefore, 
         \begin{equation}
             \notag
             {\rm div}(A^\e(x) \nabla u^{(2)}(x)_\e-A^{\hom} \nabla \umin_0(x)) = \e^2 {\rm div}(\widetilde{R}_\e(x)),
         \end{equation}
    with $\widetilde{R}_\e=(\widetilde{R}_{1,\e}, \dots, \widetilde{R}_{N,\e}).$
    Now, 
            \begin{equation}
                \notag
                \|{\rm div}(A^\e \nabla u^{(2)}_\e-A^{\hom} \nabla \umin_0) \|_{H^{-1}(\Omega)} = \|{\rm div}\widetilde{R}_\e\|_{H^{-1}(\Omega)} \leq \|\widetilde{R}_\e\|_{L^{2}(\Omega)} \leq C \e^2,
            \end{equation}
    which concludes the proof.
    \end{proof}


  \section{Compactness result in Theorem \ref{thm:mainresult}}  
   \label{sect:compactness}
  
In this section, we are going to prove a compactness result stated in Theorem \ref{thm:mainresult}(i).

    \begin{prop}
    \label{prop:cpt}    
    Let $\{\e_n\}_{n}$ be a sequence such that $\e_n\to 0$ as $n\to\infty$.
    Suppose that Assumptions {\rm (H1)}-{\rm (H3)} hold.
    If $\{u_n\}_{n}\subset H^1(\Omega)$ is a sequence such that
         \begin{equation}
             \notag
             \sup_{n} F^1_n(u_n)<\infty,
         \end{equation}
    then, $\{u_n\}_{n}$ converges to $\umin_0\in H^1(\Omega)$ weakly in $H^1(\Omega)$.
    \end{prop}

    \begin{proof}
        First, recall that thanks to Assumption {\rm (H3)}, the homogenized matrix $A_{\hom}$ satisfies the same growth conditions with the same constants, i.e.,
            \begin{align}
                \notag
                \alpha|\xi|^2\leq A_{\hom}\xi\cdot\xi\leq \beta |\xi|^2,
            \end{align}
        for all $\xi\in\R^N$.
        Let $\lambda>0$ that will be fixed later. Note that
        \begin{equation}
        \label{eq:quadratic_bound}
        -ab \geq -\frac{a^2}{2\lambda^2} - \frac{\lambda^2}{2}b^2,
        \end{equation}
        for all $a,b\geq 0$.
        Let $C_\Omega>0$ be the Poincar\'{e} constant of $\Omega$, \emph{i.e.} such that,
        \begin{equation}\label{eq:Poi}
        \|v\|_{L^2(\Omega)} \leq C_\Omega \|\nabla v\|_{L^2(\Omega;\R^N)},
        \end{equation}
        for all $v\in H^1_0(\Omega)$.
        Recall that if $u\in L^2(\Omega)$ is such that $F_{n}(u)<+\infty$, then $u\in H^1_0(\Omega)$.
        By using {\rm (H3)}, \eqref{eq:quadratic_bound} and \eqref{eq:Poi}, we get
        \begin{align*}
        F_n(u)
        &\geq \alpha\|\nabla u\|^2_{L^2(\Omega;\R^N)} 
                - \frac{\lambda^2}{2}\|u\|^2_{L^2(\Omega)}
                -\frac{1}{2\lambda^2}\|f\|^2_{L^2(\Omega)}  \\
        &\geq \left(\alpha - \frac{\lambda^2}{2} C_\Omega^2 \right)
            \|\nabla u\|^2_{L^2(\Omega;\R^N)} 
                -\frac{1}{2\lambda^2}\|f\|^2_{L^2(\Omega)}.
        \end{align*}
        Namely,
        \[
        \left(\alpha - \frac{\lambda^2}{2} C_\Omega^2 \right)
            \|\nabla u\|^2_{L^2(\Omega;\R^N)} \leq F_n(u).
        \]
        Choosing
        \[
        \lambda \in \left(0, \frac{\sqrt{2\alpha}}{C_\Omega}  \right),
        \]
        yields that
        \[
        \|\nabla u\|_{L^2(\Omega;\R^N)} \leq C ( F_n(u) + 1)
            \leq C(F^1_n(u) + 1),
        \]
        where the constant $C>0$ changes all the times.
        Since $u\in H^1_0(\Omega)$, using the Poincar\'{e} inequality again, we get that
        \[
        \|u\|_{H^1(\Omega)} < C(F^1_n(u) + 1).
        \]
        Therefore, if $\{u_n\}_{n\in\N}\subset L^2(\Omega)$ is such that
        \[
        \sup_{n\in\N} F_n^1(u_n) < +\infty,
        \]
        then, there exists a subsequence $\{u_{n_k}\}_{k\in\N}$ such that $u_{n_k}\rightharpoonup v$ weakly in $H^1(\Omega)$, for some $v\in H^1(\Omega)$.
        
        We now prove that $v=\umin_0$.
        Assume by contradiction that $v\neq \umin_0$.
        Since the minimizer $\umin_0$ is unique, it follows that
        \[
        \liminf_{k\to\infty} F_{n_k}(u_{n_k}) \geq F_\hom^0(v) > F_\hom^0(\umin_0).
        \]
        Thus,
        \[
        \lim_{k\to\infty} F^1_{n_k}(u_{n_k})
            =\lim_{k\to\infty}  \frac{F_{n_k}(u_{n_k}) - F_\hom^0(\umin_0)}{\e_{n_k}} = +\infty.
        \]
        This gives the desired contradiction.
        Since the limit is unique, we also get that the full sequence converges.
    \end{proof}

\begin{remark}
Note that, in order to get compactness, we do not need to have Assumptions {\rm{(H4)}} and {\rm(H5)} in force.
\end{remark}


\section{The liminf inequality for Theorem \ref{thm:mainresult}}
\label{sect:liminf}
In this section, we prove the lower bound of Theorem \ref{thm:mainresult} (ii).

   \begin{prop}
         \label{prop:lowerbound}
         Assume that Assumptions {\rm (H1)-(H5)} hold.
         Then, for any  sequence $\{u_n\}_{n}\subset L^2(\Omega)$ converging to $\umin_0$ with respect to $L^2(\Omega)$ it holds that  
            \begin{equation}
                \notag
                   \liminf_{n\to \infty}F_n^1(u_n) \geq  F_{\hom}^1(\umin_0).  
            \end{equation}
   \end{prop}
   
    \begin{proof}
    Without loss of generality, we can assume that
    \[
    \liminf_{n\to\infty} F^1_n(u_n)<\infty.
    \]
    Let $u^{(2)}_n$ be the function defined in \eqref{def:ue2}.
    Then, it holds that 
           \begin{align*}
              \liminf_{n\to\infty} F^{1}_n (u_n)
              &= \liminf_{n\to\infty} {F_n (u_n) - F^{0}_{\hom} (\umin_0)\over \e_n}\notag\\
              &\geq \liminf_{n\to\infty} {F_n (\umin_n) - F^{0}_{\hom} (\umin_0)  \over\e_n}\notag\\
                &\geq \liminf_{n\to\infty} {F_n (\umin_n) - F_n (u^{(2)}_n)  \over\e_n}
                + \liminf_{n\to\infty}{  F_n (u^{(2)}_n)  - F^{0}_{\hom} (\umin_0)  \over \e_n}\notag\\
                & = \liminf_{n\to\infty} I^1_n + \liminf_{n\to\infty} I_n^2
                    + \liminf_{n\to\infty} I_n^3,
            \end{align*}
    where
    \[
    I^1_n\coloneqq {  F_n (\umin_n) - F_n (u^{(2)}_n)  \over \e_n},
    \]
    \[
    I^2_n\coloneqq \frac{1}{\e_n} \left[ \int_\Omega A^{\e_n}(x) \nabla u^{(2)}_n (x) \cdot \nabla u^{(2)}_n(x) \dx
        - \int_\Omega A_{\hom} \umin_0(x)\cdot \umin_0(x) \dx  \right],
    \]
    and
    \[
    I^3_n\coloneqq \frac{1}{\e_n} \int_\Omega f(x)(u^{(2)}_n(x) - \umin_0(x)) \dx.
    \]
    We now claim that
    \[
    \lim_{n\to\infty}I^1_n =0,
    \]
    and that
    \[
    \lim_{n\to\infty}(I^2_n + I^3_n ) =F_{\hom}^1(\umin_0).
    \]
    These will give the desired result.\\


    \textbf{Step 1: limit of $I^1_n$.}
    We have that
    \begin{align*}
    I^1_n &= \frac{1}{\e_n} \left[  \int_\Omega A^{\e_n}(x) \nabla u^{(2)}_n (x)\cdot \nabla u^{(2)}_n (x)\dx  
            - \int_\Omega A^{\e_n} (x)\nabla \umin_n (x)\cdot \nabla \umin_n (x)\dx  \right] \\
    &\hspace{2cm}- \frac{1}{\e_n} \int_\Omega f (x)(u^{(2)}_n(x) - \umin_n(x)) \dx.
    \end{align*}
   Writing $u^{(2)}_n = \umin_n + (u^{(2)}_n-\umin_n)$ gives
    \begin{align}\label{eq:est_nabla}
    &\frac{1}{\e_n} \left|  \int_\Omega A^{\e_n}(x) \nabla u^{(2)}_n(x) \cdot \nabla u^{(2)}_n(x) \dx  
            - \int_\Omega A^{\e_n}(x) \nabla \umin_n(x) \cdot \nabla \umin_n(x) \dx  \right| \nonumber \\
    &\hspace{1cm}=\frac{1}{\e_n} \Bigl| \int_\Omega A^{\e_n}(x) \nabla (u^{(2)}_n-\umin_n)(x)
                                            \cdot \nabla (u^{(2)}_n-\umin_n)(x) \dx \nonumber \\
    &\hspace{3cm}+\frac{2}{\e_n} \int_\Omega A^{\e_n} (x)\nabla (u^{(2)}_n-\umin_n)(x)
                                            \cdot \nabla u^{(2)}_n(x) \dx \Bigr| \nonumber \\
    &\hspace{1cm}\leq \frac{1}{\e_n} \left\| \mathrm{div}\left( A^{\e_n} \nabla (u^{(2)}_n-\umin_n) \right)  \right\|_{H^{-1}(\Omega)}
    \left( \| u^{(2)}_n-\umin_n \|_{H^1_0(\Omega)}
     + 2\| u^{(2)}_n \|_{H^1_0(\Omega)} \right) \nonumber\\
     &\hspace{1cm}\leq \e_n\left( \| \umin_n \|_{H^1_0(\Omega)} + 3\| u^{(2)}_n \|_{H^1_0(\Omega)} \right),
    \end{align}
    where the last step follows from Proposition \ref{prop:estimate}.
    Note that from {\rm (H3)}, together with the fact that $u^{(2)}_n - \umin_n \in H^1_0(\Omega)$, it holds that
    \begin{align*}
    \alpha \| \nabla (u^{(2)}_n - \umin_n) \|^2_{L^2(\Omega)} &\leq
    \int_\Omega A^{\e_n} (x) \nabla (u^{(2)}_n - \umin_n)(x)\cdot \nabla (u^{(2)}_n - \umin_n)(x)\dx \\
    &\leq \sup_{\varphi\in H^1_0(\Omega)}\int_\Omega A^{\e_n}(x) \nabla (u^{(2)}_n - \umin_n)(x)\cdot \nabla \varphi(x)\dx \\
    &=\left\| \mathrm{div}\left( A^{\e_n} \nabla (u^{(2)}_n-\umin_n) \right)  \right\|_{H^{-1}(\Omega)} \| u^{(2)}_n - \umin_n \|_{H^1_0(\Omega)}\\
    &\leq(1+C_\Omega)\left\| \mathrm{div}\left( A^{\e_n} \nabla (u^{(2)}_n-\umin_n) \right)  \right\|_{H^{-1}(\Omega)} \| \nabla (u^{(2)}_n - \umin_n) \|_{L^2(\Omega)}.
    \end{align*}
    Therefore, calling $C_\Omega$ the Poincar\'{e} constant of $\Omega$, we obtain that
    \begin{align}\label{eq:est_f}
    \left|\frac{1}{\e_n} \int_\Omega f(x) (u^{(2)}_n(x) - \umin_n(x)) \dx \right|
    &\leq \frac{1}{\e_n} \| f\|_{L^2(\Omega)} \| u^{(2)}_n - \umin_n \|_{L^2(\Omega)} \nonumber \\
    &\leq \frac{C_\Omega}{\e_n} \| f\|_{L^2(\Omega)} \| \nabla (u^{(2)}_n - \umin_n) \|_{L^2(\Omega)}
        \nonumber \\
    &\leq \frac{(1+C_\Omega)C_\Omega}{\alpha\e_n} \| f\|_{L^2(\Omega)}
            \left\| \mathrm{div}\left( A^{\e_n} \nabla (u^{(2)}_n-\umin_n) \right)  \right\|_{H^{-1}(\Omega)} \nonumber \\
    &\leq \frac{(1+C_\Omega)C_\Omega}{\alpha} \| f\|_{L^2(\Omega)} \e_n,
    \end{align}
    where the last step follows by Proposition \ref{prop:estimate}.
    Since 
    \[
    \sup_{n} \left( \| \umin_n \|_{H^1_0(\Omega)} + 3\| u^{(2)}_n \|_{H^1_0(\Omega)} \right) < \infty,
    \]
    from \eqref{eq:est_nabla} and \eqref{eq:est_f}, we get the desired result.\\


    \textbf{Step 2: limit of $I^2_n$.}
    Note that 
    \begin{align}
        \nabla u^{(2)}_{\e_n}(x) &= \nabla \left( \umin_0(x) + {\e_n}\sum_{i=1}^{N}\psi_i^{\e_n}(x)\partial_i \umin_0(x) + {\e_n}^2\sum_{r,s=1}^{N} \chi_{rs}^{\e_n}(x)\partial^2_{rs} \umin_0(x) - {\e_n}^2\sum_{l=1}^N \psi_l^{\e_n}(x) \partial_l \widetilde{u}_1(x)  \right) \notag\\
        &= \nabla \umin_0(x) +\sum_{i=1}^{n} \nabla_y\psi_i^{\e_n}(x)\partial_i \umin_0(x)\notag\\
        &\quad +{\e_n}\left[\sum_{i=1}^{N} \psi_i^{\e_n}(x)\nabla (\partial_i \umin_0(x)) + \sum_{r,s=1}^{N} \nabla_y(\chi^{\e_n}_{rs}(x))\partial^2_{rs}\umin_0(x) - \sum_{l=1}^{N} \nabla_y(\psi^{\e_n}_l(x))\partial_l\widetilde{u}_1(x) \right] \notag \\
        &\quad+o(\e_n).\notag
    \end{align}
    Therefore,
      \begin{align}
      \notag
          \int_\Omega A^{\e_n}(x) \nabla u^{(2)}_{\e_n}(x)\cdot \nabla u^{(2)}_{\e_n}(x) \dx
            = H_n + \e_n G_n + o(\e_n),
      \end{align}
    where $H_{\e_n}$ and $G_{\e_n}$ are defined as
     \begin{align*}
         H_n&\coloneqq \int_\Omega A^{\e_n}(x) \left(\nabla \umin_0(x) +\sum_{i=1}^{N} \nabla_y(\psi_i^{\e_n}(x))\partial_i \umin_0(x) \right) \cdot \\
         &\hspace{4cm}\cdot\left(\nabla \umin_0(x) +\sum_{i=1}^{n} \nabla_y(\psi_i^{\e_n}(x))\partial_i \umin_0(x) \right)\dx,
     \end{align*}
    and
    \begin{align*}
        G_n&\coloneqq 2\int_\Omega A^{\e_n}(x)  \biggl(\nabla \umin_0(x) +\sum_{i=1}^{N} \nabla_y(\psi_i^{\e_n}(x))\partial_i \umin_0(x) \biggr) \cdot \biggl(  \sum_{i=1}^{N} \psi_i^{\e_n}(x)\nabla (\partial_i \umin_0(x))\notag\\
        &\quad \quad \quad  \quad \quad + \sum_{r,s=1}^{N} \nabla_y(\chi^{\e_n}_{rs}(x))\partial^2_{rs}\umin_0(x) 
        - \sum_{l=1}^{N} \nabla_y(\psi^{\e_n}_l(x))\partial_l\widetilde{u}_1(x)  \biggr)\dx,
    \end{align*}
    respectively.\\


    \textbf{Step 2.1: asymptotic behavior of $H_n$.}
    Using Proposition \ref{prop:RL}, and recalling that we are assuming the barycenter of $Y$ to be at the origin, we get that
    \begin{align}\label{eq:as_H_1}
    &\int_\Omega A^{\e_n}(x) \nabla \umin_0(x) \cdot \nabla \umin_0(x) \dx \nonumber \\
    &\hspace{2cm}= \sum_{i,j=1}^N \int_\Omega a_{ij}^{\e_n}(x) \partial_i \umin_0(x)\partial_j \umin_0(x)\dx 
            \nonumber \\
    &\hspace{2cm}=\sum_{i,j=1}^N \langle a_{ij}\rangle_Y\,\int_\Omega \partial_i \umin_0(x)\partial_j \umin_0(x)\dx 
                \nonumber \\
    &\hspace{3cm}+\e_n \sum_{i,j=1}^N \int_\Omega \nabla(\partial_i \umin_0\partial_j \umin_0)(x)\dx 
        \cdot \langle y a_{ij} \rangle_Y
        +o(\e_n).
    \end{align}
    Moreover,
    \begin{align}\label{eq:as_H_2}
    &\int_\Omega A^{\e_n}(x) \left( \sum_{i=1}^{N} \nabla_y(\psi_i^{\e_n}(x))\partial_i \umin_0(x) \right) \cdot \nabla \umin_0(x) \dx \nonumber \\
    &\hspace{2cm}= \sum_{i,j,s=1}^N \int_\Omega a_{js}^{\e_n}(x) \partial_{y_s} \psi_i^{\e_n}(x)
            \partial_i \umin_0(x) \partial_j \umin_0(x)\dx 
            \nonumber \\
    &\hspace{2cm}=\sum_{i,j,s=1}^N \langle a_{js} \partial_{y_s} \psi_i\rangle_Y\,
        \int_\Omega \partial_i \umin_0(x)\partial_j \umin_0(x)\dx 
                \nonumber \\
    &\hspace{3cm}+\e_n \sum_{i,j,s=1}^N \int_\Omega \nabla(\partial_i \umin_0\partial_j \umin_0)(x)\dx 
        \cdot \langle y a_{js}\partial_{y_s} \psi_i\rangle_Y
        +o(\e_n).
    \end{align}
    Finally,
    \begin{align}\label{eq:as_H_3}
    &\int_\Omega A^{\e_n}(x) \left( \sum_{i=1}^{N} \nabla_y(\psi_i^{\e_n}(x))\partial_i \umin_0(x) \right) \cdot \left( \sum_{j=1}^{N} \nabla_y(\psi_j^{\e_n}(x))\partial_j \umin_0(x) \right) \dx \nonumber \\
    &\hspace{1cm}= \sum_{i,j,r,s=1}^N \int_\Omega a_{rs}^{\e_n}(x) \partial_{y_s} \psi_i^{\e_n}(x)
            \partial_i \umin_0(x) \partial_{y_r} \psi_j^{\e_n}(x)
            \partial_j \umin_0(x)\dx 
            \nonumber \\
    &\hspace{1cm}=\sum_{i,j,r,s=1}^N \langle a_{rs}\partial_{y_s} \psi_i
        \partial_{y_r} \psi_j\rangle_Y\,
        \int_\Omega \partial_i \umin_0(x)\partial_j \umin_0(x)\dx 
                \nonumber \\
    &\hspace{2cm}+\e_n \sum_{i,j,r,s=1}^N \int_\Omega \nabla(\partial_i \umin_0\partial_j \umin_0)(x)\dx 
        \cdot \langle y a_{rs}\partial_{y_s} \psi_i\partial_{y_r} \psi_j\rangle_Y
        +o(\e_n).
    \end{align}
    Therefore, from \eqref{eq:as_H_1}, \eqref{eq:as_H_2}, \eqref{eq:as_H_3}, and the definition of $A_{\hom}$ (see \eqref{eq:def_A_hom}) we obtain that
    \begin{align}\label{eq:lim_Hn}
    H_n &= \int_\Omega A_{\hom} \nabla \umin_0(x)\cdot \nabla \umin_0(x) \dx
    \nonumber \\
    &\hspace{1cm}+ \e_n \sum_{i,j=1}^N \int_\Omega
    \nabla(\partial_i \umin_0\partial_j \umin_0)(x)\dx
    \cdot\left\langle y \left[ a_{ij}
            + 2A e_j \cdot \nabla \psi_i
            + A \nabla \psi_i \cdot \nabla \psi_j \right]\right\rangle_Y
        +o(\e_n).
    \end{align}
    \vspace{0.5cm}


    \textbf{Step 2.2: limit of $G_n$.}
Therefore,  
        \begin{align}
            \notag
            \lim_{n\to\infty} {F_{n} (\umin_n) - F^0_{\hom}(\umin_0)\over \e_n } = \liminf_{n\to\infty} G_{\e_n} (\umin_n)
        \end{align}
To make the computations more clear, we split $G_{n} (\umin_n)\coloneqq G_{n} (\umin_0, \widetilde{u}_1)$ as follows
   \begin{align*}
       G_{n} (\umin_0, \widetilde{u}_1) &\coloneqq G^1_{n} (\umin_0, \widetilde{u}_1) + G^2_{n} (\umin_0, \widetilde{u}_1) + G^3_{n} (\umin_0, \widetilde{u}_1) + G^4_{n} (\umin_0, \widetilde{u}_1) \notag\\
       &\quad + G^5_{n} (\umin_0, \widetilde{u}_1) + G^6_{n} (\umin_0, \widetilde{u}_1),
   \end{align*}
where the functionals $G^i_{n} (\umin_0, \widetilde{u}_1)$, for $i=1,\dots,6$ are given by
    \begin{align}
        \notag
        G^1_{n} (\umin_0, \widetilde{u}_1) \coloneqq 2\int_\Omega A^{\e_n}(x) \nabla \umin_0(x)\cdot \sum_{i=1}^{N}\psi_i^{\e_n}(x)\nabla (\partial_i \umin_0(x))\dx,
    \end{align}
     \begin{align}
        \notag
        G^2_{n} (\umin_0, \widetilde{u}_1)\coloneqq 2\int_\Omega A^{\e_n}(x) \nabla \umin_0(x)\cdot \sum_{r,s=1}^{N}\nabla_y(\chi_{rs}^{\e_n}(x))\partial^2_{rs} \umin_0(x)\dx,
    \end{align}
     \begin{align}
        \notag
       G^3_{n} (\umin_0, \widetilde{u}_1)\coloneqq -2\int_\Omega A^{\e_n}(x) \nabla \umin_0(x)\cdot \sum_{l=1}^{N}\nabla_y(\psi_l^{\e_n}(x))\partial_l \widetilde{u}_1(x)\dx,
    \end{align}
     \begin{align}
        \notag
        G^4_{n} (\umin_0, \widetilde{u}_1)\coloneqq 2\int_\Omega A^{\e_n}(x) \sum_{j=1}^{N} \nabla_y(\psi_j^{\e_n}(x))\partial_j \umin_0(x)\cdot\sum_{i=1}^{N}\psi_i^{\e_n}(x)\nabla (\partial_i \umin_0)\dx,
    \end{align}
     \begin{align}
        \notag
        G^5_{n} (\umin_0, \widetilde{u}_1)\coloneqq 2\int_\Omega A^{\e_n}(x) \sum_{j=1}^{N}\nabla_y(\psi_j^{\e_n}(x)) \partial_j \umin_0(x)\cdot \sum_{r,s=1}^{N}\nabla_y(\chi_{rs}^{\e_n}(x))\partial^2_{rs} \umin_0(x)\dx,
    \end{align} 
     \begin{align}
        \notag
        G^6_{n} (\umin_0, \widetilde{u}_1)\coloneqq -2\int_\Omega A^{\e_n}(x) \sum_{j=1}^{N}\nabla(\psi_j^{\e_n}(x)) \partial_j \umin_0(x)\cdot \sum_{l=1}^{N} \nabla_y(\psi_l^{\e_n}(x))\partial_l \widetilde{u}_1(x)\dx.
    \end{align}
Now, we separately compute the limit as $n\to \infty$ of each functionals $G^i_{n}$. Using the Riemann-Lebesgue lemma (see Lemma \ref{lem:RL}), we get
     \begin{align}\label{eq:G1}
         \lim_{n\to \infty} G^1_{n} (\umin_0, \widetilde{u}_1) &= 2\lim_{n\to \infty} \sum_{i=1}^N \int_\Omega \psi_i^{\e_n}(x) A^{\e_n}(x) \nabla \umin_0(x)\cdot \partial_i\nabla \umin_0(x)\dx\notag\\
         &= 2 \sum_{i=1}^N \int_\Omega  \langle \psi_i A\rangle_Y \nabla \umin_0(x)\cdot \partial_i\nabla\umin_0(x)\dx
     \end{align}
Regarding the second functional $G^2_{n}$, we obtain that 
    \begin{align}\label{eq:G2}
         \lim_{n\to \infty} G^2_{n} (\umin_0, \widetilde{u}_1) &= 2\lim_{n\to \infty} \sum_{r,s =1}^N \int_{\Omega}  \nabla \umin_0(x)\cdot A^{\e_n}(x)\nabla_y(\chi^{\e_n}_{rs}(x))\partial^2_{rs} \umin_0(x)\dx \notag\\
         &= 2\lim_{n\to \infty} \sum_{r,s =1}^N \sum_{j =1}^N \int_{\Omega}  \partial_j \umin_0(x) \left(\nabla_y(\chi^{\e_n}_{rs}(x))\cdot A^{\e_n}(x)e_j\right)\partial^2_{rs} \umin_0(x)\dx \notag\\
         &= 2\sum_{r,s =1}^N \sum_{j =1}^N \int_{\Omega} \partial_j \umin_0(x) \langle Ae_j\cdot\nabla\chi_{rs}\rangle_Y\partial^2_{rs} \umin_0(x)\dx.
    \end{align}
Here, we have exploited the symmetry of the matrix $A$ to deduce that
       \begin{align}
           \nabla \umin_0(x)\cdot A^{\e_n}(x) \nabla_y(\chi^{\e_n}_{rs}(x))
           &= \sum_{j=1}^N \partial_j\umin_0(x) \left( A^{\e_n}(x) \nabla_y\chi^{\e_n}_{rs}(x)\right)_j\notag\\
           &= \sum_{j=1}^N \partial_j\umin_0(x) \left( A^{\e_n}(x) \nabla_y\chi^{\e_n}_{rs}(x)\cdot e_j\right)\notag\\
           &= \sum_{j=1}^N \partial_j\umin_0(x) \left( \nabla_y\chi^{\e_n}_{rs}(x)\cdot  A^{\e_n}(x) e_j\right).\notag
       \end{align}
Similarly for $G^3_{n}$, it follows that
    \begin{align}\label{eq:G3}
        \lim_{n\to \infty} G^3_{n} (\umin_0, \widetilde{u}_1)
        &= -2\lim_{n\to \infty} \sum_{l=1}^N \int_{\Omega} \nabla \umin_0(x)\cdot A^{\e_n}(x) \nabla_y(\psi^{\e_n}_l(x))\partial_l\widetilde{u}_1(x)\dx \notag\\
         &=-2\lim_{n\to \infty} \sum_{l=1}^N\sum_{j=1}^N \int_{\Omega} \partial_j \umin_0(x) \cdot [A^{\e_n}(x)e_j\cdot \nabla_y(\psi^{\e_n}_l(x))] \partial_l\widetilde{u}_1(x)\dx \notag\\
        &= -2 \sum_{l=1}^N \sum_{j=1}^N \int_{\Omega}  \partial_j \umin_0(x) 
        \langle Ae_j\cdot\nabla_y\psi_l\rangle_Y \partial_l\widetilde{u}_1(x)\dx.
    \end{align}
The limit $G^4_{n} (\umin_0, \widetilde{u}_1)$ as $n\to\infty$ reads as follows
    \begin{align}\label{eq:G4}
        \lim_{n\to \infty} G^4_{n} (\umin_0, \widetilde{u}_1)
        &= 2\lim_{n\to \infty} \sum_{j=1}^N\sum_{i=1}^N \int_\Omega  \partial_j \umin_0(x)\psi_i^{\e_n}(x) A^{\e_n}(x)\nabla_y(\psi^{\e_n}_j(x))\cdot\partial_i(\nabla \umin_0(x))\dx\notag\\
        &= 2 \sum_{j=1}^N\sum_{i=1}^N \int_{\Omega}  \partial_j \umin_0(x)\langle \psi_i A\nabla_y \psi_j \rangle_Y   \cdot  \partial_i (\nabla \umin_0(x))\dx.
    \end{align}
For $G^5_{n}$, it follows that 
    \begin{align*}
        \lim_{n\to \infty} G^5_{n} (\umin_0, \widetilde{u}_1) &= 2\lim_{n\to \infty} \sum_{j=1}^N\sum_{r,s=1}^N \int_\Omega \partial_j \umin_0(x) A^{\e_n}(x) \nabla_y(\psi_j^{\e_n}(x))\cdot \nabla_y(\chi_{rs}^{\e_n}(x)) \partial_{rs}^2 \umin_0(x)\dx\notag\\
        &= 2 \sum_{j=1}^N\sum_{r,s=1}^N \int_\Omega \partial_j \umin_0(x)\langle A \nabla_y\psi_j\cdot \nabla_y\chi_{rs}\rangle_Y \partial_{rs}^2 \umin_0(x)\dx.
    \end{align*}
The variational formulation of the problem of the corrector $\psi_j$ (see \eqref{eq:first_order_corrector}) with $\chi_{rs}$ as a test function yields
          \begin{align}
             \langle A \nabla_y(\psi_j)\cdot \nabla_y(\chi_{rs}) \rangle_Y
             = \int_Y A(y)\nabla_y\psi_j(y)\cdot \nabla_y \chi_{rs}(y)\dif y = -\int_Y A(y)e_j\cdot \nabla_y\chi_{rs}(y)\dif y
             = -\langle A e_j\cdot \nabla_y\chi_{rs} \rangle_Y.\notag
          \end{align}
This implies that
    \begin{align}\label{eq:G5}
      \lim_{n\to \infty} G^5_{n} (\umin_0, \widetilde{u}_1)
      =- 2 \sum_{j=1}^N\sum_{r,s=1}^N \int_\Omega \partial_j \umin_0(x)\langle A e_j\cdot \nabla_y\chi_{rs} \rangle_Y \partial_{rs}^2 \umin_0(x)\dx.
    \end{align}  
Finally, the limit of of $G^6_{n}$ as $n \to \infty$ is
    \begin{align}
        \lim_{n\to \infty} G^6_{n} (\umin_0, \widetilde{u}_1) &= -2\lim_{n\to \infty} \sum_{j=1}^N\sum_{l=1}^N \int_{\Omega} \partial_j \umin_0(x) A^{\e_n}(x) \nabla_y(\psi_j^{\e_n}(x))\cdot \nabla_y(\psi_l^{\e_n}(x))\partial_l \widetilde{u}_1(x)\dx\notag\\
        &=-2 \sum_{j=1}^N\sum_{l=1}^N \int_{\Omega} \partial_j \umin_0(x) \langle A \nabla_y\psi_j\cdot \nabla_y\psi_l \rangle_Y\partial_l \widetilde{u}_1(x)\dx.\notag
    \end{align}
Using again the variational formulation of the problem of the corrector $\psi_j$ (see \eqref{eq:first_order_corrector}), choosing as test function $\psi_l$, yields to 
          \begin{align}
             \langle A \nabla_y\psi_j\cdot \nabla_y\psi_l \rangle_Y 
             &= -\langle A e_j\cdot \nabla_y\psi_l \rangle_Y.\notag
          \end{align}
Thanks to this equality, it follows that 
           \begin{align}\label{eq:G6}
           \lim_{n\to \infty} G^6_{n} (\umin_0, \widetilde{u}_1) = 2 \sum_{j=1}^N\sum_{l=1}^N \int_{\Omega} \partial_j \umin_0(x) \langle A e_j\cdot \nabla_y\psi_l \rangle_Y\partial_l \widetilde{u}_1(x)\dx.
          \end{align}    
Gathering formulas \eqref{eq:G1}-\eqref{eq:G6} and noting that \eqref{eq:G2} and \eqref{eq:G5} as well as \eqref{eq:G3} and \eqref{eq:G6} cancel out, we deduce that  
    \begin{align}\label{eq:lim_Gn}
        &\lim_{n\to \infty} G_{n} (\umin_0, \widetilde{u}_1) = \lim_{n\to \infty} G^1_{n} (\umin_0, \widetilde{u}_1) + \lim_{n\to \infty} G^4_{n} (\umin_0, \widetilde{u}_1)\notag\\
         &\quad= 2 \sum_{i=1}^N \int_\Omega   \langle \psi_i A\rangle_Y \nabla \umin_0(x)\cdot \partial_i(\umin_0(x))\dx  + 2 \sum_{j=1}^N\sum_{i=1}^N \int_{\Omega} \partial_j \umin_0(x) \langle \psi_i A\nabla_y\psi_j \rangle_Y  \cdot \partial_i (\nabla \umin_0(x)) \dx.
    \end{align}


    \textbf{Step 3: limit of $I^3_n$.}
    A direct computation shows that 
            \begin{align}
            \label{eq:lim_f}
                \lim_{n\to\infty} {1\over \e_n}\int_\Omega f(x)(u^{(2)}_n(x) - u_0^{\min}(x))\dx &= \lim_{n\to\infty} \sum_{j=1}^N\int_\Omega f(x) \psi_j^{\e_n}(x) \partial_ju_0^{\min}(x)\dx \notag\\
                &= \sum_{j=1}^N \langle \psi_j \rangle_Y \int_\Omega f(x) \partial_j u_0^{\min}(x)\dx. 
            \end{align}
    \end{proof}


\section{The limsup inequality for Theorem \ref{thm:mainresult}}
\label{sect:limsup}
In this section, we prove the upper bound of Theorem \ref{thm:mainresult} (ii).

    \begin{prop}
        \label{prop:upperbound}
        Assume Assumptions {\rm (H1)-(H5)} to hold.
        Then, there exists a sequence $\{u_n\}_{n}\subset H^1_0(\Omega)$ converging to $\umin_0$ weakly in $H^1(\Omega)$ such that
        \[
            \lim_{n\to\infty} F^1_n(u_n) = F^1_{\hom}(\umin_0).
        \]
    \end{prop}

    \begin{proof}
        We define, for $n\in\N$,
        \[
        u_n(x)\coloneqq \umin_0(x) + \e\sum_{i=1}^{N}\psi_i^{\e_n}(x)\partial_i \umin_0(x).
        \]
        The weak convergence of $\{u_n\}_n$ to $\umin_0$ is a direct computation.
        Moreover, the desired convergence of the energy follows from \eqref{eq:G1}, \eqref{eq:G3}, \eqref{eq:G4}, and \eqref{eq:lim_f}.
    \end{proof}


\section{Proof of Theorem \ref{thm:main_Lp}}
\label{sect:generalfunc}

This section is devoted to the proof of Theorem \ref{thm:main_Lp}.
We first illustrate the idea of the proof by considering the special case
\[
V(x,p)=a(x)|p|^2,
\]
where $a:\R^N\to\R$ is an $Y$-periodic function such that $0<c_1\leq a(c)\leq c_2<+\infty$ for all $x\in Y$.
Then, it holds that
\[
V_{\hom}(p) = \min\left\{ \int_Y a(y)|p+\varphi(y)|^2\dy \,:\, \varphi\in L^2(Y),\, \int_Y\varphi(y)\dy=0 \right\}.
\]
This mininization problem admits a solution, that we can compute explicitly. Indeed, the Euler-Lagrange equation gives the existence of a constant $c(p)\in\R$ such that
\[
2a(y)(p+\varphi(y)) = c(p),
\]
for almost every $y\in Y$. Recalling that $a(y)>0$, we obtain
\[
\varphi(y) = \frac{c(p)}{2a(y)} - p.
\]
Moreover, the constant $c(p)$ can be computed by the zero average requirement for $\varphi$:
\[
c(p) = 2p \left( \int_Y \frac{1}{a(y)}\dy \right)^{-1}.
\]
Thus,
\[
V_{\hom}(p) = p^2 \left( \int_Y \frac{1}{a(y)}\dy \right)^{-1},
\]
which yields
\[
G_\hom(v) = \int_\Omega v^2(x)\dx\, \left( \int_Y \frac{1}{a(y)}\dy \right)^{-1},
\]
for all $v\in L^2(\Omega)$.

Fix $m\in\R$. We now consider the solution $v^{\min}_0\in L^2(\Omega)$ to the minimization problem
\[
\min\left\{ G_\hom(v) \,:\, v\in L^2(\Omega),\, \int_\Omega v(x)\dx = m  \right\}.
\]
The Euler-Lagrange equation gives the existence of a constant $c\in\R$ such that
\[
2v^{\min}_0(x) \left( \int_Y \frac{1}{a(y)}\dy \right)^{-1} = c,
\]
for almost every $x\in\Omega$. Thus, $v^{\min}_0$ is constant. The mass constraint gives that
\[
v^{\min}_0(x) = \frac{m}{|\Omega|},
\]
for all $x\in\Omega$, which yields
\begin{equation}\label{eq:G_hom_min}
G_\hom(v^{\min}_0) = {m^2\over |\Omega|}\left( \int_Y \frac{1}{a(y)}\dy \right)^{-1}.\\
\end{equation}

In a similar way, we can obtain the solution $v^{\min}_n$ to the minimization problem
\[
\min\left\{ G_{\varepsilon_n}(v) \,:\, v\in L^2(\Omega),\, \int_\Omega v(x)\dx = m  \right\}.
\]
This gives
\[
v^{\min}_n(x) = \frac{c_n}{2a^{\varepsilon_n}\left(x\right)},
\]
where
\[
c_n\coloneqq 2m\left( \int_\Omega \frac{1}{a^{\varepsilon_n}\left(x\right)}\dx \right)^{-1}.
\]
Thus,
\begin{equation}\label{eq:G_n_min}
G_n(v^{\min}_n) = m^2 \left( \int_\Omega \frac{1}{a^{\varepsilon_n}\left(x\right)}\dx \right)^{-1}.
\end{equation}
Note that thanks to our assumptions on $\Omega$ and $\varepsilon_n$, we get
\[
\int_\Omega \frac{1}{a^{\varepsilon_n}\left(x\right)}\dx
	= |\Omega| \int_Y \frac{1}{a(y)} \dy.
\]
Therefore, from \eqref{eq:G_hom_min} and \eqref{eq:G_n_min}, we conclude that
\[
G_n(v^{\min}_n) = G_{\hom}(v^{\min}_0),
\]
for all $n\in\N\setminus\{0\}$, as desired.\\

In the general case, we use a similar argument as the one implemented above. 
We first consider the homogenized energy density
\[
V_{\hom}(p) = \min\left\{ \int_Y V(y,p+\varphi(y))\dy \,:\, \varphi\in L^2(Y),\, \int_Y\varphi(y)\dy=0 \right\}.
\]
Since we are assuming $p\mapsto V(y,p)$ to be strictly convex, for each $y\in Y$ we denote by $\partial_p V^{-1}(y):\R\to\R$ the inverse of the map $\partial_p V(y, \cdot)$, and we use the notation $v\mapsto \partial_p^{-1}V(y)[v]$ to avoid the use of too many round parenthesis.
We get that the optimal perturbation $\varphi$ satisfies
\[
\varphi(y) = \partial_p^{-1}V(y)[c(p)] - p,
\]
for some $c(p)\in\R$, which can be computed by using the zero average constrained
\[
p=\int_Y \partial_p^{-1}V(y)[c(p)]\dy.
\]
If we now consider the solution $v^{\min}_0\in L^2(\Omega)$ to the minimization problem
\[
\min\left\{ G_\hom(v) \,:\, v\in L^2(\Omega),\, \int_\Omega v(x)\dx = m  \right\},
\]
we get that
\[
v^{\min}_0 = \frac{m}{|\Omega|}.
\]
In particular,
\[
G_{\hom}(v^{\min}_0) = |\Omega|\int_Y V\left(y,\partial^{-1}_pV(y)\left[c\left( \frac{m}{|\Omega|} \right)\right]\right) \dy.\\
\]

If we now consider the solution $v^{\min}_n$ to the minimization problem
\[
\min\left\{ G_{n}(v) \,:\, v\in L^2(\Omega),\, \int_\Omega v(x)\dx = m  \right\},
\]
we get that
\[
\partial_pV\left( \frac{x}{\varepsilon_n}, v^{\min}_n(x) \right)=c_n,
\]
where, using the mass constraint,
\[
m = \int_\Omega \partial_p^{-1}V\left(\frac{x}{\varepsilon_n}\right)[c_n]\dx
=|\Omega| \int_Y \partial_p^{-1}V(y)[c_n]\dy.
\]
This gives that
\[
G_{n}(v^{\min}_n)
	= \int_\Omega V\left( \frac{x}{\varepsilon_n}, \partial^{-1}_pV\left(\frac{x}{\varepsilon_n}\right)[c_n] \right)\dx
	=|\Omega| \int_Y V\left( y, \partial^{-1}_pV\left(y\right)[c_n] \right)\dy.
\]
It is possible the choose
\[
c_n = \frac{m}{|\Omega|},
\]
which gives $G_{n}(v^{\min}_n) = G_{\hom}(v^{\min}_0)$ for all $n\in\N\setminus\{0\}$ as desired.


\subsection*{Acknowledgment}
We would like to thank George Allaire for useful conversation about the subject.
Moreover, the authors would like to thank CIRM Luminy for its hospitality during the Research in Residence.
RC was partially supported under NWO-OCENW.M.21.336, MATHEMAMI - Mathematical Analysis of phase
Transitions in HEterogeneous MAterials with Materials Inclusions. LD acknowledges support of the Austrian Science Fund (FWF)
projects 10.55776/ESP1887024 and 10.55776/Y1292.


\end{document}